\documentclass[preprint,12pt]{elsarticle}

\usepackage{amsfonts}
\usepackage{amsmath}
\usepackage{amssymb}
\usepackage{mathtools}
\usepackage{graphicx}
\usepackage{epstopdf}
\usepackage[ruled,vlined,algo2e]{algorithm2e}
\usepackage{bbm}
\usepackage{lscape}
\usepackage{yfonts}
\usepackage{verbatim}
\usepackage{xcolor}
\usepackage{hyperref}
\usepackage{xurl}
\usepackage{subcaption}
\usepackage{algorithm}
\usepackage{algorithmic}
\usepackage{amsthm}

\ifpdf
  \DeclareGraphicsExtensions{.eps,.pdf,.png,.jpg}
\else
  \DeclareGraphicsExtensions{.eps}
\fi

\DeclareSymbolFont{lettersA}{U}{txmia}{m}{it}
\DeclareMathSymbol{\R}{\mathord}{lettersA}{"92}
\DeclareMathSymbol{\C}{\mathord}{lettersA}{"83}
\newcommand{\ffast}{\ensuremath{f^{\mathfrak{f}}}}
\newcommand{\fmed}{\ensuremath{f^{\mathfrak{m}}}}
\newcommand{\fslow}{\ensuremath{f^{\mathfrak{s}}}}
\newcommand{\fimpl}{\ensuremath{f^{\mathfrak{i}}}}
\newcommand{\fexpl}{\ensuremath{f^{\mathfrak{e}}}}
\newcommand{\mfs}{\mathfrak{s}}
\newcommand{\mff}{\mathfrak{f}}

\newcommand{\bigO}[1]{\ensuremath{\mathop{}\mathopen{}\mathcal{O}\mathopen{}\left(#1\right)}}

\newcommand{\Hh}{\ensuremath{H\text{-}h}}
\newcommand{\HTol}{\ensuremath{H\text{-}Tol}}
\newcommand{\Decoupled}{\text{\emph{Decoupled}}}
\definecolor{mygreen}{rgb}{0.13, 0.55, 0.13}
\allowdisplaybreaks

\newcommand{\corrA}[1]{{\color{black}#1}}
\newcommand{\corrB}[1]{{\color{black}#1}}

\newproof{pf}{Proof}
\newdefinition{remark}{Remark}
\newdefinition{definition}{Definition}
\newdefinition{assumption}{Assumption}
\newtheorem{theorem}{Theorem}[section]
\newtheorem{lemma}[theorem]{Lemma}

\makeatletter
\newcommand*{\addFileDependency}[1]{
  \typeout{(#1)}
  \@addtofilelist{#1}
  \IfFileExists{#1}{}{\typeout{No file #1.}}
}
\makeatother

\ifpdf
\hypersetup{
  pdftitle={Efficient and Flexible Multirate Temporal Adaptivity},
  pdfauthor={D.~R. Reynolds, S. Amihere, D. Mitchell, and V.~T. Luan}
}
\fi

\journal{Journal of Computational and Applied Mathematics}

\begin{document}

\begin{frontmatter}

\title{Efficient and Flexible Multirate Temporal Adaptivity\tnoteref{label1}}

\affiliation[affil1]{organization={Department of Mathematics and Statistics, University of Maryland Baltimore County},
    city={Baltimore},
    state={Maryland},
    country={USA}
}

\affiliation[affil2]{organization={Department of Mathematics and Statistics, Texas Tech University},
    city={Lubbock},
    state={Texas},
    country={USA}
}

\tnotetext[label1]{
This work was funded in part by the U.S. Department of Energy, Office of Science, Office of Advanced Scientific Computing Research, Scientific Discovery through Advanced Computing (SciDAC) Program through the FASTMath Institute, under DOE award DE-SC0021354.
V.T. Luan was partially supported by NSF awards DMS–2309821 and DMS–2531805.
}

\cortext[cor1]{Corresponding author}

\author[affil1]{Daniel R. Reynolds\corref{cor1}}
\ead{dreynolds@umbc.edu}

\author[affil1]{Sylvia Amihere}
\ead{samihere@umbc.edu}

\author[affil1]{Dashon Mitchell}
\ead{dashonm1@umbc.edu}

\author[affil2]{Vu Thai Luan}
\ead{Vu.Luan@ttu.edu}

\begin{abstract}
In this work we present two new families of multirate time step adaptivity controllers, that are designed to work with embedded multirate infinitesimal (MRI) time integration methods for adapting time steps when solving problems with multiple time scales.  We compare these controllers against competing approaches on two benchmark problems, showing that the proposed methods offer dramatically improved performance and flexibility.  The combination of embedded MRI methods and the proposed controllers enable adaptive simulations of problems with a potentially arbitrary number of time scales, achieving high accuracy while maintaining low computational cost.  Additionally, we introduce a new set of embeddings for the family of explicit multirate exponential Runge--Kutta (MERK) methods of orders 2 through 5, resulting in the first-ever fifth-order embedded MRI method.  Finally, we compare the performance of a wide range of embedded MRI methods on our benchmark problems to provide guidance on how to select an appropriate MRI method and multirate controller.
\end{abstract}

\begin{keyword}
  multirate time integration \sep time step adaptivity \sep initial-value problems
  \MSC[2010]  65L05 \sep 65L06 \sep 65L20
\end{keyword}

\end{frontmatter}


\section{Introduction}
\label{sec:introduction}

In this work, we focus on time-step adaptivity within multirate infinitesimal (MRI) methods that are used when solving systems of ordinary differential equation initial value problems \corrB{(IVPs)} of the form
\begin{equation}
\label{eq:multirateIVP}
\begin{split}
    y'(t) &= \fslow(t,y) + \ffast(t,y), \quad t>t_0\\
    y(t_0) &= y_0,
\end{split}
\end{equation}
Here, the operator $\ffast$ contains physical processes that evolve on a fast time scale with typical time step size $h$, and $\fslow$ contains processes that evolve on a slower time scale with typical step size $H \gg h$.  Multirate methods are frequently used for problems of this form when $\fslow$ is considerably more costly to evaluate than $\ffast$, and thus algorithms that evaluate $\fslow$ infrequently have the potential for significant runtime improvements over single rate approaches that evolve all processes with time step $h$.

Generally speaking, an explicit infinitesimal method for evolving a single time step $y_n\approx y(t_n)$ to $y_{n+1}\approx y(t_n+H_n)$, with its embedded solution $\tilde{y}_{n+1}\approx y(t_n+H_n)$, proceeds through the following sequence \cite{schlegel_multirate_2009,sanduClassMultirateInfinitesimal2019,chinomonaImplicitExplicitMultirateInfinitesimal2021,luanNewClassHighOrder2020,luanMultirateExponentialRosenbrock2022,fishImplicitExplicitMultirate2024}.
\begin{enumerate}
    \item Let: $Y_1=y_n$
    \item For $i=2,...,s$:
      \label{eq:mristep_stages}
      \begin{enumerate}
      \item Solve: $v_i'(\theta) = \ffast(\theta,v_i) + r_i(\theta)$, for $\theta\in[\theta_{0,i},\theta_{F,i}]$ with $v_i(\theta_{0,i})=v_{0,i}$.\label{eq:mristep_faststage}
      \item Let: $Y_i = v_i(\theta_{F,i})$.
      \end{enumerate}
    \item Solve: $\tilde{v}_s'(\theta) = \ffast(\theta,\tilde{v}_s) + \tilde{r}_{s}(\theta)$, for $\theta\in[\theta_{0,s},\theta_{F,s}]$ with $\tilde{v}_s(\theta_{0,s})=v_{0,s}$.\label{eq:mristep_fastembedding}
    \item Let: $y_{n+1} = Y_s$ and $\tilde{y}_{n+1} = \tilde{v}_s(\theta_{F,s})$.
\end{enumerate}
MRI methods are determined by their fast stage time intervals $[\theta_{0,i},\theta_{F,i}]$, initial conditions $v_{0,i}$, forcing functions $r_i(\theta)$, and embedding forcing function $\tilde{r}_{s}(\theta)$. Both $r_i(\theta)$ and $\tilde{r}_{s}(\theta)$ are constructed using linear combinations of $\left\{\fslow(\theta_{F,j},Y_j)\right\}$, and serve to propagate information from the slow to the fast time scales.  Implicit and implicit-explicit extensions of the above MRI algorithm exist, that involve replacing step \ref{eq:mristep_faststage} with an implicit solve for some internal stages $Y_i$. The embedded solution $\tilde{y}_{n+1}$ is similar to embedded solutions in standard Runge--Kutta methods, in that it approximates the solution to an alternate order of accuracy and may be computed with minimal extra effort beyond computing $y_{n+1}$; in this case it is computed by re-solving the last stage with an alternate forcing function, $\tilde{r}_s(\theta)$.

As seen in steps \ref{eq:mristep_faststage} and \ref{eq:mristep_fastembedding} above, computation of each stage in an explicit MRI method requires solving a secondary, \emph{inner}, IVP.
Typically, these inner IVPs are not solved exactly, and instead are approximated using another numerical method with inner time steps $h_{n,m}$ such that $\max_m |h_{n,m}| \ll |H_n|$.  \corrB{Our goal in this paper is to examine optimal ``controller'' algorithms that can adaptively select both $H_n$ and $h_{n,m}$ within MRI methods.}

In the sections that follow, we first summarize how temporal adaptivity is performed for problems with a single time scale.  In Section \ref{sec:multirate_controllers} we introduce previous work, and propose two new approaches, for dynamically adapting both the slow and fast time steps, $H_n$ and $h_{n,m}$, within MRI methods.  Two of these multirate controller families require estimates of the accumulated fast temporal error, so we describe our algorithms for computing those estimates in Section \ref{sec:fast_error_estimate}.  The majority of this manuscript focuses on numerical tests of these multirate controllers in Section \ref{sec:numerical_tests}.  In that section, we first describe our two benchmark problems that we use throughout our numerical tests.  We then outline the set of embedded MRI methods that we will test, including new embeddings for the explicit \textit{multirate exponential Runge--Kutta} (MERK) methods that result in the first ever fifth-order embedded MRI method.  The remainder of Section \ref{sec:numerical_tests} performs a variety of numerical experiments to examine the proposed multirate adaptivity strategies and the embedded MRI methods.  We discuss our conclusions from these tests and identify future directions to extend this work in Section \ref{sec:conclusions}.

\section{Single Rate Control}
\label{sec:singlerate_control}

Before introducing multirate temporal adaptivity, we briefly review the typical approaches for single-rate adaptivity, that attempt to control the \emph{local error} induced within a given time step.  \corrB{For additional information on this topic, we refer the reader to the references \cite{hairer_solving_1993,hindmarsh_sundials_2005,gardner_enabling_2022}.}  Assuming that the initial condition at the time step $t_n$ is exact, i.e., $y_n - y(t_n) = 0$, the local error is defined to be
\begin{equation}
    \label{eq:local-error}
    \varepsilon_n = y(t_n+h_n)-y_{n+1}.
\end{equation}
Single-rate controllers then adapt the step size $h_n$ to achieve two primary objectives: (a) ensure the local error is ``small enough'', i.e.,
\begin{equation}
    \label{eq:local-error-test}
    \| \varepsilon_n \| < 1,
\end{equation}
where $\|\cdot\|$ is a norm that incorporates the requested temporal accuracy\corrB{, such as the \emph{weighted root-mean-square} (WRMS) norm from \cite{hindmarsh_sundials_2005,gardner_enabling_2022},
\begin{equation}
  \label{eq:WRMS_norm}
  \| \varepsilon_n \| = \left(\frac{1}{N} \sum_{j=1}^N \left(\frac{|\varepsilon_{n,j}|}{\text{reltol} |y_{n,j}| + \text{abstol}_j}\right)^2\right)^{1/2}
\end{equation}
for IVPs with $N$ components, where $\text{reltol}$ and $\text{abstol}_j$ are user-specified relative and absolute tolerances, respectively,} and (b) perform as few time steps as possible to meet objective (a).  Additional considerations include a smooth step size progression  \cite{gustafsson_control_1991,soderlind_automatic_2002,soderlind_digital_2003,soderlind_time-step_2006}, and improving the efficiency and robustness of iterative algebraic solvers that may be used at each step \cite{gustafsson_control-theoretic_1994,gustafsson_control_1997}.

Single rate adaptivity may be easily incorporated into standard ``time marching'' schemes.  Instead of utilizing a fixed step size $h_n = h$, one begins at $t_0$ with an initial step size $h_0$, and initializes the counter $n=0$.  Then in a loop over step attempts:
\begin{itemize}
\item[(i)] compute a candidate approximation $y_{n+1}^*\approx y(t_n+h_n)$ and an approximation of the local error $\varepsilon_n^* = \tilde{y}_{n+1}-y_{n+1}^* \approx y(t_n+h_n)-y_{n+1}^*$.
\item[(ii)] If $\|\varepsilon_n^*\|<1$ set $t_{n+1} = t_n+h_n$, $y_{n+1}=y_{n+1}^*$, and $n = n+1$.
Else: reject the approximation $y_{n+1}^*$ (and the corresponding step size $h_n$).
\item[(iii)] Use $\|\varepsilon_n^*\|$ to estimate a step size $\tilde{h}$ for the subsequent step attempt.
\end{itemize}

Whether $\varepsilon_n^*$ passes or fails the error test \eqref{eq:local-error-test}, step (iii) must select $\tilde{h}$ to use on the next step attempt.  This \emph{adaptivity controller} is a function that typically depends on a small set of $(h_{n-k},\varepsilon_{n-k}^*)$ pairs \corrB{and the asymptotic order of accuracy for the method $p$}, i.e.,
\[
   \tilde{h} = \mathcal{C}(h_n,\varepsilon_n^*,h_{n-1},\varepsilon_{n-1}^*,\ldots,p).
\]
There are a myriad of choices for controllers $\mathcal{C}$ in the literature, but the simplest approach is the so-called \emph{I controller}.  \corrB{Since $y_{n+1}$ has order $p$, then from \cite[Sec.~II.3]{hairer_solving_1993} there exists a function $c(t)$ independent of $h_n$ such that},
\begin{align}
  \notag
  &\varepsilon_n \approx y(t_n+h_n)-y_{n+1}^* = c(t_n)\,h_n^{p+1} + \bigO{h_n^{p+2}} = \bigO{h_n^{p+1}}\\
  \notag
  \Rightarrow\quad&\\
  \label{eq:asymptotic-error}
  &\|\varepsilon_n\| \approx \|c(t_n)\|\,h_n^{p+1}.
\end{align}
Since we have just computed the estimate $\|\varepsilon_n^*\|$ based on a step with size $h_n$, a piecewise constant approximation $c(t)\approx c_n$ allows us to ``solve'' for $\|c_n\| = \|\varepsilon_n^*\| / h_n^{p+1}$.
Assuming that $c(t)$ will remain relatively constant from one step attempt to the next, the candidate step size $\tilde{h}$ is computed to exactly attain an error goal of 1 (or any $tol$):
\begin{align*}
  tol = \|c_n\|\tilde{h}^{p+1}
  \quad\Rightarrow\quad
  \tilde{h} = h_n\left(\frac{tol}{\|\varepsilon_n^*\|}\right)^{1/(p+1)}.
\end{align*}
This formula automatically grows or shrinks the step size based on whether the error test $\|\varepsilon_n^*\|^{1/(p+1)}\le tol$ passed or failed, respectively.  Due to the multiple assumptions above, a safety factor $\sigma<1$ is typically applied,
\[
  \tilde{h} = \sigma h_n\left(\frac{tol}{\|\varepsilon_n^*\|}\right)^{1/(p+1)}.
\]
More advanced approaches for $\mathcal{C}$ are typically based on control theory, and use additional $(h_{n-k},\varepsilon_{n-k}^*)$ values to build higher-degree piecewise polynomial approximations of the principal error function  \cite{gustafsson_control_1991,gustafsson_control-theoretic_1994,gustafsson_control_1997,soderlind_automatic_2002,soderlind_digital_2003,soderlind_time-step_2006}.

To extend these ideas to the multi-rate context, we require two complementary components: (a) strategies to estimate local temporal error, and (b) algorithms for selecting step sizes following each attempted step.  Since (b) dictates our needs for (a), we begin with multirate temporal control.

\section{Multirate Temporal Controllers}
\label{sec:multirate_controllers}

Multirate temporal adaptivity has not been widely studied in the literature.  The two articles we are aware of are \cite{sarshar_design_2019} and \cite{fishAdaptiveTimeStep2023}. The former is specific to ``MrGARK'' methods \cite{gunther_multirate_2016}, that require a more rigid algorithmic structure than the MRI algorithm introduced in Section \ref{sec:introduction}.  The latter works with MRI methods, and will serve as a baseline comparison for our work.  In the following subsections, we consider three such approaches.

\subsection{Coupled Stepsize (``\Hh'') Control}
\label{sec:Hh-control}

The adaptive multirate controllers from \cite{fishAdaptiveTimeStep2023} simultaneously predict both a slow step size $H_n$, and a multirate ratio $M_n$, such that the inner solver takes fixed small substeps $h_n = H_n/M_n$ throughout each subinterval $[\theta_{0,i},\theta_{F,i}]$.  Due to their use of fixed inner steps that are adapted only when the outer step $H_n$ is adapted, we will refer to these as \emph{coupled stepsize (``\Hh'') controllers}.  In our subsequent comparisons, we consider four controllers introduced in \cite{fishAdaptiveTimeStep2023}: MRI-CC, MRI-LL, MRI-PI, and MRI-PID.

\corrA{We note, however, that after introducing these controllers in \cite{fishAdaptiveTimeStep2023}, we found that they struggled on applications where the optimal ``multirate ratio'' $M_n = H_n / h_n$ was large.  Specifically, we found that when $M_n \le 10$ these controllers performed well (consistent with the test problems in \cite{fishAdaptiveTimeStep2023}), but their performance deteriorated for $M_n \ge 20$.  In the following subsection, we briefly explain this phenomenon.}

\subsubsection{\corrA{Asymptotic analysis of \Hh\ controller limitations}}
\label{subsubsec:Hh-analysis}

\corrA{As an exemplar of the \Hh\ controller family, we analyze the \emph{Constant-Constant H-h} Controller (MRI-CC) introduced in \cite{fishAdaptiveTimeStep2023} in the limit of large multirate ratio.  We note that MRI-CC mimics the \emph{I} controller from Section \ref{sec:singlerate_control}, but updates $M_n$ using a formula that couples the fast and slow error estimates from the most recent time step,
\begin{equation}
  \label{eq:MRICC}
  \begin{split}
    H^{CC}_{n+1} & = H_{n} \|\varepsilon^{\mfs}_{n}\|^{-\frac{k_1}{P}}\\
    M^{CC}_{n+1} & = M_{n} \|\varepsilon^{\mfs}_{n}\|^{-\frac{(p+1)k_{1}}{Pp}} \|\varepsilon^{\mff}_{n}\|^{\frac{k_{2}}{p}},
  \end{split}
\end{equation}
where $k_1$ and $k_2$ are controller parameters, $P$ and $p$ are the local orders of accuracy of the embeddings for the slow and fast solvers, and $\varepsilon^{\mfs}_{n}$ and $\varepsilon^{\mff}_{n}$ the estimated local errors for the time step $t_n \to t_{n}+H_n$, respectively.

By comparison, if separate \emph{I} controllers were performed independently at each time scale, the corresponding multirate ratio would be
\begin{equation}
  \label{eq:MRI_I}
  M^I_{n+1} = M_{n}\|\varepsilon^{\mfs}_{n}\|^{-\frac{1}{P}} \|\varepsilon^{\mff}_{n}\|^{\frac{1}{p}}.
\end{equation}

To understand the limiting case, it suffices to consider a problem with decoupled fast and slow time scales with $\|\varepsilon^{\mfs}_{n}\| = 1$ and $\|\varepsilon^{\mff}_{n}\| \gg 1$, which should result in $H_{n+1} = H_n$ and $M_{n+1} \gg M_n$.  Taking the limit as $\|\varepsilon^{\mff}_{n}\|\to \infty$ of the ratio $M_{n+1}^{CC} / M_{n+1}^{I}$, we have
\begin{align*}
  \lim_{\|\varepsilon_{n}^{\mff}\| \to \infty} \frac{M_{n+1}^{CC}}{M_{n+1}^{I}}
  &= \|\varepsilon_{n}^{\mfs}\|^{-\frac{1}{P} \left( \frac{(p+1)k_{1}}{p} - 1 \right)} \lim_{\|\varepsilon_{n}^{\mff}\| \to \infty} \|\varepsilon_{n}^{\mff}\|^{\frac{k_{2}-1}{p}}.
\end{align*}

The parameters $k_{1}=0.42$ and $k_{2}=0.44$ are determined in \cite{fishAdaptiveTimeStep2023} by optimizing controller performance across a suite of test problems.  Inserting these values and noting that $p>0$, we see that
\begin{equation*}
  \lim_{\|\varepsilon_{n}^{\mff}\| \to \infty} \frac{M_{n+1}^{CC}}{M_{n+1}^{I}}
  = \|\varepsilon_{n}^{\mfs}\|^{-\frac{1}{P} \left( \frac{0.42(p+1)}{p} - 1 \right)} \lim_{\|\varepsilon_{n}^{\mff}\| \to \infty} \|\varepsilon_{n}^{\mff}\|^{-\frac{0.56}{p}} = 0,
\end{equation*}
indicating that although the physical problem may necessitate large $M_{n+1}$ values (and thus very small fast steps $h_n$), the MRI-CC controller will artificially limit this value, potentially leading to computational inefficiency or inaccuracy.  This effect will be readily visible in the numerical results in Section \ref{sec:numerical_tests} below.}

\subsection{``\Decoupled'' Multirate Control}
\label{sec:decoupled-control}

Our simplest proposed MRI controller uses two decoupled single-rate controllers to separately adapt the inner and outer time steps, i.e.,
\begin{equation}
    \label{eq:decoupled-control}
    \begin{split}
      \tilde{H} &= \mathcal{C}^{\mfs}(H_n,\varepsilon_n^{\mfs},H_{n-1},\varepsilon_{n-1}^{\mfs},\ldots,P),\\
      \tilde{h} &= \mathcal{C}^{\mff}(h_{n,m},\varepsilon_{n,m}^{\mff},h_{n,m-1},\varepsilon_{n,m-1}^{\mff},\ldots,p),
    \end{split}
\end{equation}
where $P$ and $p$ are the orders of the MRI method and inner solver, respectively.  Here, $(H_k,\varepsilon_k^{\mfs})$ are the stepsize and local error estimates for time step $k$ at the slow time scale, and $(h_{k,\ell},\varepsilon_{k,\ell}^{\mff})$ are the stepsize and local error estimates for the fast substep $\ell$ within the slow step $k$.  We note that in this approach, both $\mathcal{C}^{\mfs}$ and $\mathcal{C}^{\mff}$ are decoupled, and thus selection of $\tilde{H}$ and $\tilde{h}$ occurs independently through application of any pair of single-rate controllers.  We summarize their use in the following pseudocode.

Given the current state $y_n$, candidate step $H_n$, and controllers $\mathcal{C}^{\mfs}$ and $\mathcal{C}^{\mff}$, perform a single MRI step attempt as follows.
\begin{enumerate}
  \item Let: $Y_1 = y_n$.
  \item For each MRI stage $i = 2,\ldots,s$:
    \begin{enumerate}
    \item Use an adaptive solver with $\mathcal{C}^{\mff}$ to ensure $\|\varepsilon_{n,m}^{\mff}\| \le 1$ for the IVP $v_i'(\theta) = \ffast(\theta,v_i) + r_i(\theta)$, $\theta\in[\theta_{0,i},\theta_{F,i}]$, $v_i(\theta_{0,i})=v_{0,i}$.
    \item Let $Y_i = v_i(\theta_{F,i})$.
    \end{enumerate}
  \item Use an adaptive solver with $\mathcal{C}^{\mff}$ to ensure $\|\varepsilon_{n,m}^{\mff}\| \le 1$ for the IVP $\tilde{v}_s'(\theta) = \ffast(\theta,\tilde{v}_s) + \tilde{r}_s(\theta)$, $\theta\in[\theta_{0,s},\theta_{F,s}]$, $\tilde{v}_s(\theta_{0,s})=v_{0,s}$.
  \item Let: $y_{n+1}^* = Y_s$, $\tilde{y}_{n+1} = \tilde{v}_s(\theta_{F,s})$, and $\varepsilon_{n}^{\mfs} = \tilde{y}_{n+1}-y^*_{n+1}$.
  \item Use $\mathcal{C}^{\mfs}$ to select a new step size $\tilde{H}$ for the ensuing step attempt.
\end{enumerate}

Since this approach ignores any coupling between time scales, we expect these controllers to work well for multirate problems wherein the time scales are somewhat decoupled, where errors introduced at one scale do not ``pollute'' the other, for example reaction-diffusion systems \cite{Ropp2004}, acoustic-elastic wave systems \cite{Petersson2018}, or transport-chemistry models \cite{Mcrae1982}.  Furthermore, due to its decoupled nature, this technique can be trivially extended to an arbitrary number of time scales, allowing temporal adaptivity for so-called ``telescopic'' multirate methods.

\subsection{Stepsize-tolerance (``\HTol'') Control}
\label{sec:HTol-control}

Our second proposed family are the so-called ``\HTol'' multirate controllers.  As outlined in Section \ref{sec:singlerate_control}, standard controllers attempt to adapt step sizes to control \emph{local error} within each step to achieve a requested tolerance.  However, MRI methods must ask another ``inner'' adaptive fast-scale solver to produce the stage solutions $v_i(\theta_{F,i})$, that result from \corrB{numerically integrating} over each interval $[\theta_{0,i},\theta_{F,i}]$.  Local errors within the fast integrator may accumulate, resulting in an overall fast-solver error $\varepsilon^{\mff}_{i,true}$ that could exceed the requested tolerance.  If that inner solver can produce \emph{both} $v_i(\theta_{F,i})$ and an estimation of the accumulated error, $\varepsilon^{\mff}_{i}$, then the tolerances provided to that inner solver can be adjusted to ensure MRI stage solutions $Y_i = v_i(\theta_{F,i})$ that are within the overall tolerances requested of the outer MRI method.

To this end, first assume that at step $n$, \corrB{we define $\text{tolfac}_n$ as a multiplicative factor indicating how much more accurate the inner solver should be than the MRI method}, i.e., if the MRI method strives to achieve local errors $\|\varepsilon_n^{\mfs}\| \le 1$, then the inner solver strives to achieve local substep errors $\|\varepsilon_{n,m}^{\mff}\| \le \text{tolfac}_n$.  \corrB{Additionally assume that the inner solver is run in ``local extrapolation'' mode \cite[Sec.~II.7]{hairer_solving_1993}, i.e., the solution is one order more accurate than the embedding and thus $\|\varepsilon_{n,m,true}^{\mff}\| \approx h_{n,m}\|\varepsilon_{n,m}^{\mff}\|$.}  \corrB{Then accounting for accumulation of temporal error in the inner solver, the overall fast error over the slow step $[t_n, t_n+H_n]$ has the form
\begin{equation}
   \label{eq:fast_error_accumulation_assumption}
   \|\varepsilon^{\mff}_{n}\| \le \chi_n H_n \text{tolfac}_n,
\end{equation}
where $\chi_n$ is independent of $\text{tolfac}_n$, but may vary in time.  We prove this result in Lemma \ref{lem:accumulated_fast_error}, after defining appropriate notation.  At stage $i$ of slow step $n$, the inner
solver approximates solutions to
\begin{equation}\label{eq:fast_ivp_n}
  v_i'(\theta)= f^{\mff}(\theta, v_i(\theta)) + r_i(\theta),\quad  v_i(\theta_{0,i})=v_{0,i},\quad
  \theta\in[\theta_{0,i},\theta_{F,i}],
\end{equation}
where the forcing function $r_i(\theta)=r(\theta; t_i,y_i)$ is constructed from slow-stage data at preceding stages $j<i$ in the step.  Let $v_{i,m}$ denote the numerical approximation produced by the adaptive inner solver at the inner nodes
$\theta_m$ (with $\theta_0=\theta_{0,i}$ and $\theta_M=\theta_{F,i}$), i.e.,
$v_{i,m}\approx v_i(\theta_m)$ for
\[
  \theta_{0,i}=\theta_0<\theta_1<\cdots<\theta_m<\cdots<\theta_M=\theta_{F,i},
\]
with substep sizes $h_m=\theta_{m+1}-\theta_m$. We assume the inner solver is a one-step method, so that
\begin{equation}\label{eq:one_step_Phi}
  v_{i,m+1}= v_{i,m} + h_m \Psi^{\mff}_{h_m}(\theta_m,v_{i,m})=: \Phi_{h_m}(\theta_m,v_{i,m}),\quad m=0,\dots,M-1,
\end{equation}
for some one-step method $\Phi_h$ applied to \eqref{eq:fast_ivp_n}.
Here, $\Psi^{\mff}_{h_m}(\theta_m,v_{i,m})$ approximates $f^{\mff}(\theta, v_n(\theta)) + r_i(\theta)$ (possibly using
multiple internal stages).\\
Denote the global and local errors after $m$ steps respectively as
\begin{equation}\label{eq:global_local_errors}
  e^{\mff}_{i,m} := v_{i,m} - v_i(\theta_{m}), \quad
  \varepsilon_{i,m}^{\mff} := \Phi_{h_m}\big(\theta_m,v_i(\theta_m)\big) - v_i(\theta_{m+1}).
\end{equation}
\begin{lemma}
\label{lem:accumulated_fast_error}
Assume that $f^{\mff}(\theta,v)$ is locally Lipschitz continuous in $v$, uniformly in $\theta$ with constant
$L_f$ and that the inner error controller satisfies the local bound
\begin{equation}\label{eq:error_control_n}
  \|\varepsilon_{i,m,true}^{\mff}\| \le  h_m \text{\em tolfac}_n, \qquad m=0,\dots,M-1.
\end{equation}
Then the accumulated fast error at the end of the inner solve \eqref{eq:one_step_Phi} satisfies
\begin{equation}\label{eq:accum_error_bound_n}
  \|e^{\mff}_{i,m}\|  \le  \chi_i\,(\theta_{F,i}-\theta_{0,i})\,\text{\em tolfac}_n,
\end{equation}
where the constant $\chi_i$ is independent of $\text{\em tolfac}_n$.\\
In particular, if the MRI method satisfies $\sum_{i=1}^s \left(\theta_{F,i}-\theta_{0,i}\right) = c H_n$ for a fixed constant $c$, then the accumulated error over the $s$ fast IVP solves satisfies
\begin{equation}\label{eq:accum_error_bound_Hn}
  \|\varepsilon_n^{\mff}\|\le \chi_n H_n\text{\em tolfac}_n,
\end{equation}
where $\chi_n = c \max_i \chi_i$.
\end{lemma}
\begin{proof}
First, we observe that
\begin{align*}
  e^{\mff}_{i,m+1}
  &= v_{i,m+1} - v_i(\theta_{m+1}) \\
  & = \Phi_{h_m}(\theta_m,v_{i,m}) -  \Phi_{h_m}\big(\theta_m,v_i(\theta_m)\big)  + \Phi_{h_m}\big(\theta_m,v_i(\theta_m)\big)  -  v_i(\theta_{m+1})\\
  & = \Phi_{h_m}(\theta_m,v_{i,m}) - \Phi_{h_{m}}(\theta_{m},v_i(\theta_{m}))  + \varepsilon_{i,m}^{\mff}.
\end{align*}
Under the assumption that $f^{\mff}(\theta,v)$ has the  Lipschitz constant $L_{f}$, it follows from \eqref{eq:one_step_Phi} that  $\Phi_h(\theta,v)$ also satisfies the Lipschitz condition (in its second argument) with the Lipschitz constant $L_{\Phi_{h_m}} = 1 + h_m L_f$. Therefore,
\[
  \|e^{\mff}_{i,m+1}\| \le (1+L_{f}h_{m})\|e^{\mff}_{i,m}\|+\| \varepsilon_{i,m}^{\mff}\|.
\]
Solving this recursion, applying a discrete Gr\"onwall inequality, and using the fact that  $\sum_{m=0}^{M-1}h_{m}=\theta_{F,i}-\theta_{0,i}$, yields
\begin{align*}
  \|e^{\mff}_{i,M}\| \le \exp\Bigl(L_{f}\sum_{m=0}^{M-1}h_{m}\Bigr)\sum_{m=0}^{M-1}\| \varepsilon_{i,m}^{\mff}\|
   = \exp\bigl(L_{f}(\theta_{F,i}-\theta_{0,i})\bigr))\sum_{m=0}^{M-1}\| \varepsilon_{i,m}^{\mff}\|.
\end{align*}
Using \eqref{eq:error_control_n} gives
\[
  \|e^{\mff}_{i,M}\| \le \exp\bigl(L_{f}(\theta_{F,i}-\theta_{0,i})\bigr)  \sum_{m=0}^{M-1} h_m \text{tolfac}_n
  \
  =  \chi_i\,(\theta_{F,i}-\theta_{0,i})\,\text{tolfac}_n
\]
with $\chi_i = \exp\bigl(L_{f}(\theta_{F,i}-\theta_{0,i})\bigr)$.
This proves \eqref{eq:accum_error_bound_n}, and consequently, \eqref{eq:accum_error_bound_Hn}.
\end{proof}
}


\corrB{To convert equation \eqref{eq:fast_error_accumulation_assumption} into a multirate controller, we replace the inequality with an equality, and then regroup terms} as $\|\varepsilon^{\mff}_{n}\| = \left(\chi_n H_n\right) \left(\text{tolfac}_n\right)^{0+1}$, we see that this fits the standard asymptotic error assumption used in single-rate controllers \eqref{eq:asymptotic-error}, through identifying the control parameter $\text{tolfac}_n$ with $h_n$, $\chi_n H_n$ with $\|c(t_n)\|$ \corrB{(note that this is independent of $\text{tolfac}_n$)}, and by defining the ``order'' $p=0$.  Thus, any single-rate controller that relies on the asymptotic error assumption \eqref{eq:fast_error_accumulation_assumption} may be used to adjust $\text{tolfac}_n$ between slow step attempts.

Thus we construct an \HTol\ controller from three single-rate controllers:
\begin{itemize}
\item $\mathcal{C}^{\mfs,H}$ -- adapts $H_n$ to achieve user-requested solution tolerances using data $(H_n, \varepsilon^{\mfs}_n, H_{n-1}, \varepsilon^{\mfs}_{n-1}, \ldots, P)$.
\item $\mathcal{C}^{\mfs,Tol}$ -- adapts $\text{tolfac}_n$ using the strategy described above with data $(\text{tolfac}_n, \varepsilon^{\mff}_n, \text{tolfac}_{n-1}, \varepsilon^{\mff}_{n-1}, \ldots, 0)$.
\item $\mathcal{C}^{\mff}$ -- adapts inner time steps $h_{n,m}$ to achieve the current requested tolerance, $\text{tolfac}_n$, using data $(h_{n,m},\varepsilon^{\mff}_{n,m},h_{n,m-1},\varepsilon^{\mff}_{n,m-1},\ldots,p)$.
\end{itemize}
We summarize their use in the following pseudocode.

Given the current state $y_n$, candidate step $H_n$, and controllers $\mathcal{C}^{\mfs,H}$, $\mathcal{C}^{\mfs,Tol}$ and $\mathcal{C}^{\mff}$, perform a single MRI step attempt as follows.
\begin{enumerate}
  \item Let: $Y_1 = y_n$.
  \item For each MRI stage $i = 2,\ldots,s$:
    \begin{enumerate}
    \item Use an adaptive solver with $\mathcal{C}^{\mff}$ to ensure $\|\varepsilon_{n,m}^{\mff}\| \le \text{tolfac}_n$ for the IVP $v_i'(\theta) = \ffast(\theta,v_i) + r_i(\theta)$, $\theta\in[\theta_{0,i},\theta_{F,i}]$, $v_i(\theta_{0,i})=v_{0,i}$.
    \item Let $Y_i = v_i(\theta_{F,i})$.
    \end{enumerate}
  \item Use an adaptive solver with $\mathcal{C}^{\mff}$ to ensure $\|\varepsilon_{n,m}^{\mff}\| \le \text{tolfac}_n$ for the IVP $\tilde{v}_s'(\theta) = \ffast(\theta,\tilde{v}_s) + \tilde{r}_s(\theta)$, $\theta\in[\theta_{0,s},\theta_{F,s}]$, $\tilde{v}_s(\theta_{0,s})=v_{0,s}$.
  \item Let: $y^*_{n+1} = Y_s$, $\tilde{y}_{n+1} = \tilde{v}_s(\theta_{F,s})$, $\varepsilon_{n}^{\mfs} = \tilde{y}_{n+1} - y^*_{n+1}$, and retrieve the accumulated fast error $\varepsilon^{\mff}_n$ from the inner solver.
  \item Use $\mathcal{C}^{\mfs,H}$ to select a new step size $\tilde{H}$ for the ensuing step attempt.
  \item Use $\mathcal{C}^{\mfs,Tol}$ to select a tolerance factor $\widetilde{\text{tolfac}}$ for the ensuing step attempt.
\end{enumerate}

This class of controllers may also be applied to telescopic MRI methods, since it only focuses on the relationship between two successive time scales.

\section{Fast Temporal Error Estimation}
\label{sec:fast_error_estimate}

The \Hh\ and \HTol\ MRI adaptivity families require estimates of the temporal errors at \emph{both} the slow and fast time scales, $\varepsilon^{\mfs}_n$ and $\varepsilon^{\mff}_n$, respectively.  While the slow error may be estimated using the MRI method's time step solution and embedding, $\varepsilon^{\mfs}_n = \|y_n - \tilde{y}_n\|$, non-intrusive approaches for estimating $\varepsilon^{\mff}_n$ are more challenging.  We employ the following two strategies.

Since the \HTol\ multirate controller \corrB{relies on} an adaptive fast solver, we leverage the fact that at each substep of the fast integrator $t_{n,m}$, it computes an estimate of the local temporal error, $\varepsilon^{\mff}_{n,m}$ (e.g., using its own time step solution and embedding).  Since the fast solver may be run with a different relative tolerance than the slow solver, we denote this as $\text{reltol}^{\mff}_n$.  \corrB{As shown in \eqref{lem:accumulated_fast_error}, typical theory expects these local errors to grow exponentially based on the local Lipschitz constant, $L_n$; however in practice this accumulation is generally more well-behaved.  In this paper, we choose a moderate accumulation assumption that the total fast error may be approximated as} the sum of the substep errors,
\begin{equation}
  \label{eq:additive_accumulation}
  \varepsilon^{\mff,add}_n = \frac{\text{reltol}^{\mff}_n}{\text{reltol}^{\mfs}_n} \sum_{m\in M} \|\varepsilon^{\mff}_{n,m}\|,
\end{equation}
where $M$ is the set of all steps since the accumulator was last reset (generally at the start of the step, $t_n$).  To understand the tolerance ratio, we \corrB{assume that} both the fast and slow solvers use a weighted root-mean squared norm \corrB{\eqref{eq:WRMS_norm}}, such that \corrB{$\|\varepsilon^{\mff}_{n,m} \| \le 1$ is ``sufficiently accurate.''}  The prefactor $\text{reltol}^{\mff}_n$ serves to convert the accumulated fast errors to a raw estimate of the relative error, before the factor $1/\text{reltol}^{\mfs}_n$ scales this back to normalized relative error.

\corrB{We note that other accumulation strategies may be more well-suited to specific problems, where local errors may decay or even grow exponentially.  In exploring options, we also considered a ``maximum accumulation'' strategy,
\begin{equation}
  \label{eq:maximum_accumulation}
  \varepsilon^{\mff,max}_n = \frac{\text{reltol}^{\mff}_n}{\text{reltol}^{\mfs}_n} \max_{m\in M} \|\varepsilon^{\mff}_{n,m}\|,
\end{equation}
and an ``average accumulation'' strategy,
\begin{equation}
  \label{eq:average_accumulation}
  \varepsilon^{\mff,avg}_n = \frac{\text{reltol}^{\mff}_n}{\text{reltol}^{\mfs}_n} \sum_{m\in M} \frac{h_{n,m}}{T_M} \|\varepsilon^{\mff}_{n,m}\|,
\end{equation}
where $T_M$ is the elapsed simulation time over the set $M$.  In testing these approaches, we found that the additive strategy \eqref{eq:additive_accumulation} showed the most robust performance; however, both of the strategies \eqref{eq:maximum_accumulation} and \eqref{eq:average_accumulation} are also included in the ARKODE time integration solver \cite{reynolds_arkode_2023} from SUNDIALS \cite{hindmarsh_sundials_2005,gardner_enabling_2022}, where we have implemented our methods in this paper.}

Since the \Hh\ multirate controller uses fixed fast substeps we cannot assume that the fast solver is adaptive\corrB{, and thus it may not compute substep estimates $\|\varepsilon_{n,m}^{\mff}\|$}.  Thus, we use a more costly approach \corrB{based on Richardson extrapolation to estimate the accumulated fast error over each slow MRI stage $i$.  This} runs the fast integrator using fixed time steps of size $h$ and $kh$ \corrB{for $k\ne 1$,} (e.g., $k=2$) to achieve fast solutions \corrB{$v_{i,h}(\theta_{F,i})$ and $v_{i,kh}(\theta_{F,i})$}, which we then subtract to estimate the fast temporal error \corrB{over that slow stage},
\corrB{\[
  \varepsilon^{\mff,dbl}_i = \left|\frac{1}{k^{p}-1}\right| \; \|v_{i,h}(\theta_{F,i}) - v_{i,kh}(\theta_{F,i})\|,
\]}
where $p$ is the global order of accuracy for the \corrB{inner IVP solver \cite[Sec.~II.4]{hairer_solving_1993}.  We then take the maximum of these to estimate the overall fast temporal error,
\begin{equation}
  \label{eq:double_step_accumulation}
  \varepsilon^{\mff,dbl}_n = \max_{1\le i\le s} \varepsilon^{\mff,dbl}_i.
\end{equation}
We note that this estimate may struggle for stiff problems, where order reduction may render the effective order to differ from $p$.}

\section{Numerical Tests}
\label{sec:numerical_tests}

In this section, we perform numerical tests to verify the proposed \Decoupled\ and \HTol\ multirate adaptivity strategies, and to compare their performance against the previous \Hh\ controllers.  All codes for performing these tests and reproducing the figures in this paper are available in \cite{mri_paper_github_2025}, and use the ARKODE time integration solver \cite{reynolds_arkode_2023}, which is part of the SUNDIALS library \cite{hindmarsh_sundials_2005,gardner_enabling_2022}.

\begin{figure}[t]
    \centering
    \includegraphics[width=0.5\linewidth]{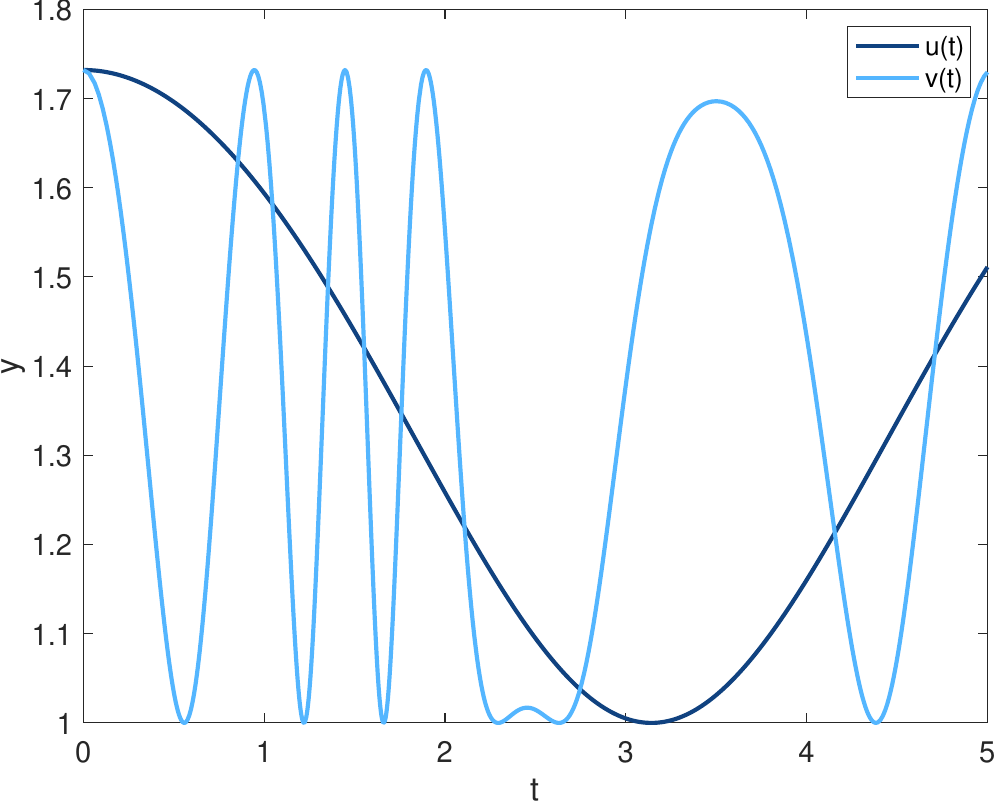}\hfill
    \includegraphics[width=0.5\linewidth]{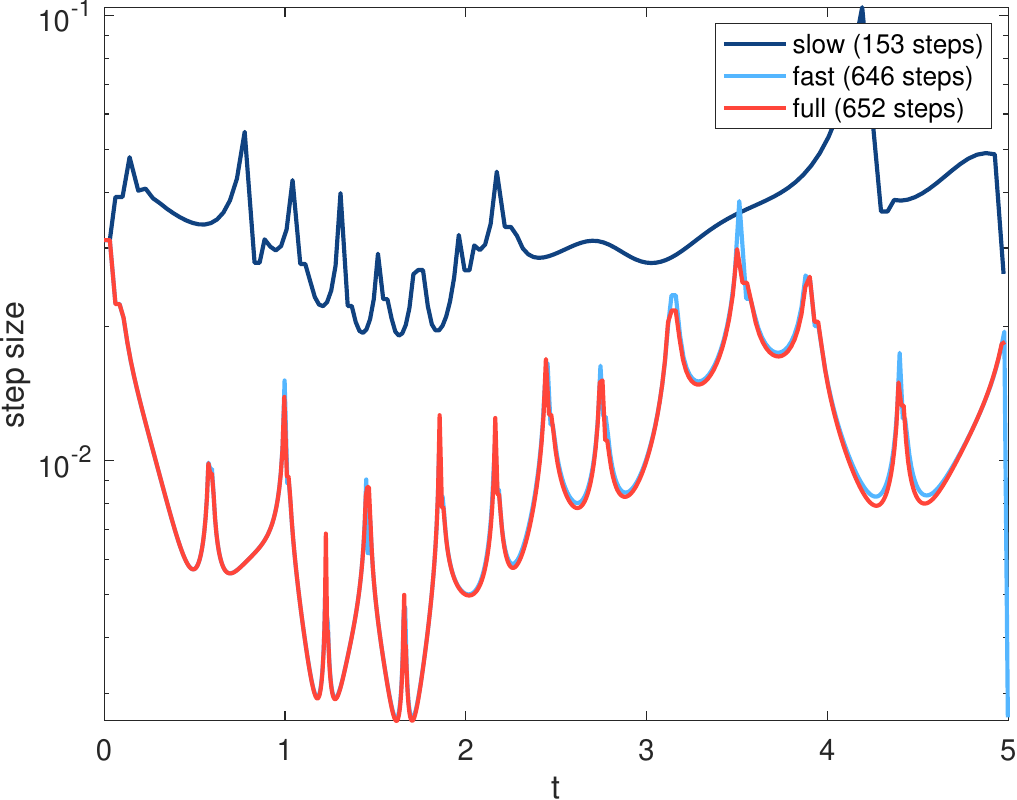}
    \caption{\corrA{Left: example solution to the 2-scale KPR problem \eqref{eq:KPR} with scale separation factor $\omega=5$.  The oscillations in $v$ change frequency through the time interval, evolving approximately $\omega$ times faster than $u$.  Right: internal time step sizes used when solving single-scale problems and full problem, showing the need for smaller step sizes when including the faster time scale.}}
    \label{fig:KPR-2scale}
\end{figure}

\corrA{In our experience, MRI methods are most often used for two types of applications.  The first are problems with multiple dynamcially-evolving time scales that must be resolved accurately, and where the fast time scale evolves considerably more rapidly than the slow time scale \cite{nonaka2026magnex}.  The second are problems with a stiff component, and where instead of applying an implicit or implicit-explicit (ImEx) solver to that term, it is subcycled with an explicit method \cite{yao_advancing_2025}.  We note that both \cite{nonaka2026magnex} and \cite{yao_advancing_2025} use MRI methods from ARKODE, but use fixed slow step sizes $H_n$ that could be improved based on the adaptive methods presented here.

To explore the performance of our adaptive methods, we test them on two benchmark problems that exhibit the aforementioned characteristics.}
The first is a multirate, nonlinear version of the Kv{\ae}rno-Prothero-Robinson (KPR) ODE test problem, that has been slightly modified from \cite[Section 6.2]{fishImplicitExplicitMultirate2024}:
\begin{equation}
  \label{eq:KPR}
  \begin{pmatrix} u'(t) \\ v'(t) \end{pmatrix}
  = \begin{bmatrix} G & e_s \\ e_f & -1 \end{bmatrix}
    \begin{pmatrix} \left(u^2-p-2\right)/(2u) \\ \left(v^2-q-2\right)/(2v) \end{pmatrix}
    + \begin{pmatrix} p'(t) / (2u) \\ q'(t) / (2v) \end{pmatrix}
\end{equation}
over $0< t < 5$, where $p(t) = \cos(t)$ and $q(t) = \cos(\omega t(1 + e^{-(t-2)^2}))$.  This problem has analytical solution $u(t) = \sqrt{2+p(t)}$ and $v(t) = \sqrt{2+q(t)}$, and its behavior is dictated by the parameters: $e_s$ determines the strength of coupling from the fast to the slow scale, $e_f$ determines the coupling strength from the slow to the fast scale, $G<0$ determines the stiffness at the slow time scale, and $w$ determines the time-scale separation factor.
\corrA{
An illustration of the analytical solutions for this problem using parameters $G=-10$, $e_s=e_f=1/10$, and a small time scale separation factor of $\omega=5$, is given in Figure \ref{fig:KPR-2scale} (left).}
We split the right-hand side above into up to three functions for the slow-implicit ($\fimpl$), slow-explicit ($\fexpl$) and fast ($\ffast$) operators, respectively,
\begin{equation}
  \label{eq:KPR-split}
  \underbrace{\begin{pmatrix} G\frac{u^2-p-2}{2u} + e_s\frac{v^2-q-2}{2v} \\ 0 \end{pmatrix}}_{\fimpl}
  + \underbrace{\begin{pmatrix} \frac{p'(t)}{2u} \\ 0 \end{pmatrix}}_{\fexpl}
  + \underbrace{\begin{pmatrix} 0 \\ e_f\frac{u^2-p-2}{2u} - \frac{v^2-q-2-q'(t)}{2v} \end{pmatrix}}_{\ffast};
\end{equation}
for non-ImEx MRI methods the combined slow operator is $\fslow = \fimpl+\fexpl$.
\corrA{
We verify that the splitting \eqref{eq:KPR-split} results in the desired time scales by running single scale simulations with $\fslow$ and $\ffast$ separately, with initial conditions $u(0)=v(0)=\sqrt{3}$ (the analytical solution values at $t=0$).  We adaptively integrate each problem using MATLAB's \texttt{ode23} method, with absolute tolerance $10^{-11}$ and relative tolerance $10^{-6}$.  In Figure \ref{fig:KPR-2scale} (right), we overlay the internal time steps taken to solve both problems, along with the internal steps required to solve the full problem, where we see that as expected these steps must adapt throughout each run, and also that the fast and slow scales are generally separated by a factor of $\omega$.}

\begin{figure}[t]
    \centering
    \includegraphics[width=0.465\linewidth]{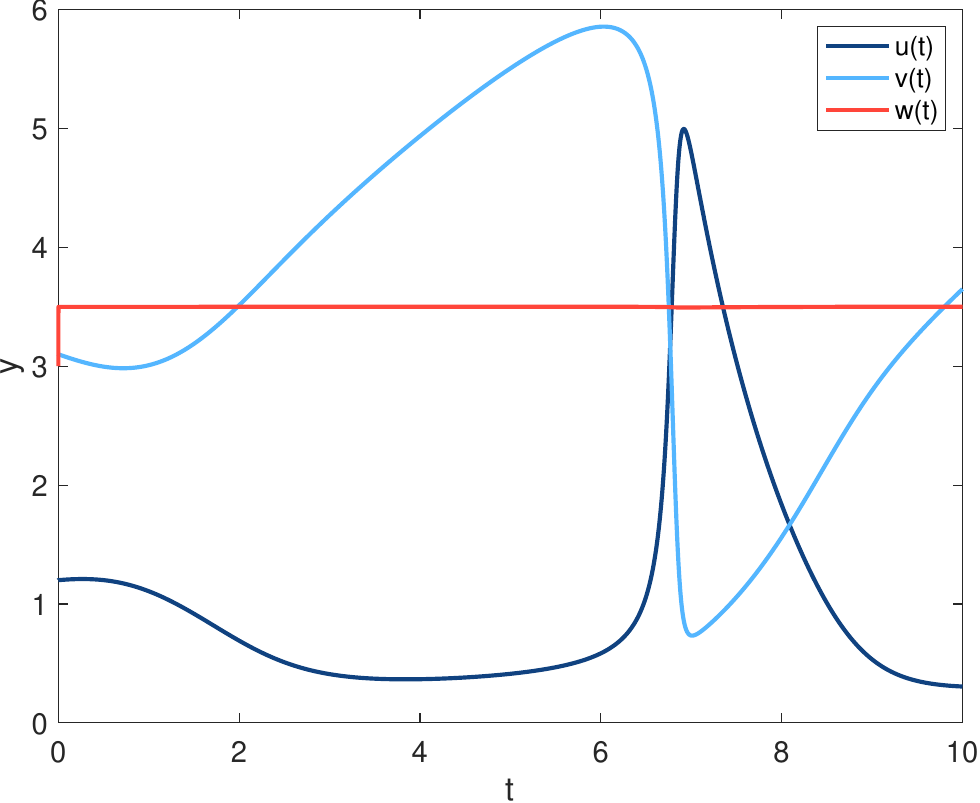}\hfill
    \includegraphics[width=0.535\linewidth]{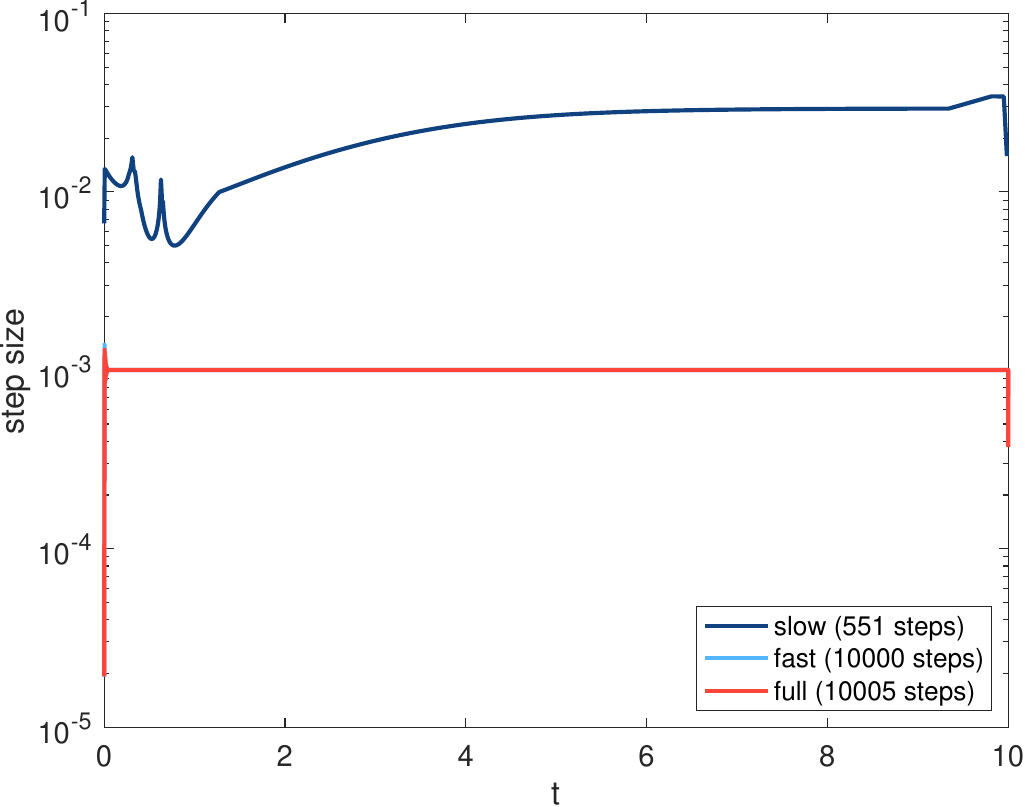}
    \caption{\corrA{Left: example solution to the stiff Brusselator problem \eqref{eq:Bruss} with stiffness factor $\epsilon=1/2500$.  Right: internal time step sizes used when solving single-scale problems and full problem.  Since the fast time steps are primarily stability-limited, they remain relatively constant, while the slow time steps must rapidly adapt at first, and then remain steady for the remainder of the time interval.}}
    \label{fig:Bruss}
\end{figure}

The second benchmark is a stiff version of the Brusselator problem, originally proposed in \cite{prigogine_symmetry_1968},
\begin{equation}
  \label{eq:Bruss}
  \begin{pmatrix} u'(t) \\ v'(t) \\ w'(t) \end{pmatrix}
  = \begin{pmatrix} a + v u^2 - (w+1)u\\ wu - v u^2\\ \frac{b-w}{\epsilon} - w u \end{pmatrix}
\end{equation}
for $0 < t < 10$. We use the initial conditions $u(0)=1.2$, $v(0)=3.1$, and $w(0) = 3$, and parameters $a=1$, $b=3.5$ and $\epsilon = 5\times 10^{-6}$, unless stated otherwise.  We note that $\epsilon$ determines the stiffness and/or time scale separation factor for the problem, such that typically $H/h = 1/(100\epsilon)$.  We again define a splitting for this right-hand side,
\begin{equation}
  \label{eq:Bruss-split}
  \underbrace{\begin{pmatrix} - (w+1)u\\ w u\\ -w u \end{pmatrix}}_{\fimpl}
  + \underbrace{\begin{pmatrix} a + v u^2\\ -v u^2\\ 0 \end{pmatrix}}_{\fexpl}
  + \underbrace{\begin{pmatrix} 0\\ 0\\ \frac{b-w}{\epsilon} \end{pmatrix}}_{\ffast}
\end{equation}
\corrA{
Figure \ref{fig:Bruss} provides an illustration of the solutions for this problem using $\epsilon=1/2500$, as well as the time steps taken when solving the slow-only, fast-only, and full problems, using an identical approach as performed for Figure \ref{fig:KPR-2scale}.}

\subsection{Embedded MRI Methods}
\label{sec:embedded-MRI}

To explore multirate temporal adaptivity, we rely on MRI methods that include embeddings, to provide approximations of both the time step solution $y_n$ and embedding $\tilde{y}_n$.  To this end, we run tests using the following 15 methods (with orders of accuracy in parentheses):
\begin{itemize}
    \item MRI-GARK methods from \cite{sanduClassMultirateInfinitesimal2019}:
    \begin{itemize}
        \item Explicit \texttt{ERK22a} (2), \texttt{ERK22b} (2), \texttt{ERK33a} (3), \texttt{ERK45a} (4);
        \item Implicit \texttt{IRK21a} (2), \texttt{ESDIRK34a} (3), and \texttt{ESDIRK46a} (4);
    \end{itemize}
    \item the explicit \texttt{RALSTON2} (2) MRI-GARK method from \cite{roberts_fast_2022};
    \item IMEX-MRI-SR methods from \cite{fishImplicitExplicitMultirate2024}: \texttt{IMEXSR21} (2), \texttt{IMEXSR32} (3), and \texttt{IMEXSR43} (4);
    \item explicit MERK methods \texttt{MERK21} (2), \texttt{MERK32} (3), \texttt{MERK43} (4), and \texttt{MERK54} (5).
\end{itemize}
The embeddings for each of the above methods have an order of accuracy one lower than the method itself.  We note that the original \texttt{MERK2}, \texttt{MERK3}, \texttt{MERK4}, and \texttt{MERK5} methods from \cite{luanNewClassHighOrder2020} do not include embeddings.  We provide embeddings for each in  \ref{sec:embeddedMERK}.

\subsection{Multirate Temporal Controllers}
\label{sec:verification_MRI_adaptivity}

With our fast temporal estimation strategies in place, we now examine the performance of the MRI adaptivity algorithms from Section \ref{sec:multirate_controllers}.  For both the \Decoupled\ and \HTol\ controller approaches, we construct MRI controllers using each of the $H211$, $H_0211$, $H_0321$, $H312$ and $I$ single-rate adaptivity controllers \cite{soderlind_digital_2003}, from the ARKODE library \cite{reynolds_arkode_2023}; we additionally test the previously-introduced \Hh\ controllers MRI-CC, MRI-LL, MRI-PI, and MRI-PID \cite{fishAdaptiveTimeStep2023}.  For all tests, we pair each of the 15 MRI methods from Section \ref{sec:embedded-MRI} with a fast explicit Runge--Kutta method of the same order.  We apply these to two test problems:
\begin{itemize}
    \item the KPR problem \eqref{eq:KPR} (with parameters $G=-100$, $e_s=5$, $e_f=0.5$, $\omega=\{50, 500\}$), to assess controller performance when both fast and slow scales evolve dynamically, and where the scale separation factor is varied;
    \item the stiff Brusselator problem \eqref{eq:Bruss} (with parameter $\epsilon=\{10^{-4}, 10^{-5}\})$, to assess controller performance when the fast scale is stability-limited at varied stiffness levels but generally evolves at a constant rate, but the slow scale varies in time.
\end{itemize}
We run each problem over their full temporal duration using a fixed absolute tolerance $\text{abstol}= 10^{-11}$, and a variety of relative tolerances, $\text{reltol}=10^{-k}, k=3,\ldots,7$.  This corresponds with a total of 4200 test combinations.  For each test combination and time step $(t_{n-1},y_{n-1}) \to (t_{n},y_{n})$, we compute a reference solution using a fifth-order single-rate explicit solver with initial condition $(t_{n-1},y_{n-1})$ and tolerances $\text{reltol}=10^{-10}$ and $\text{abstol}=10^{-12}$, respectively, to evolve to $(t_n,y_{ref}(t_n))$.  We then determine each method's ability to achieve the target local accuracy as
\begin{equation}
  \label{eq:accuracy}
  \text{accuracy} = \max \left| \frac{y_{n,l} - y_{ref,l}(t_n)}{\text{abstol} + \text{reltol}\, |y_{ref,l}(t_n)|} \right|,
\end{equation}
where the maximum is taken over all solution components ($l$) and over all time steps ($n$).  We note that an optimal method would achieve ``accuracy'' values exactly equal to 1; however in practice most combinations of adaptivity controllers and time integration methods would be considered successful if they achieve an accuracy value within a factor of 100 of the requested relative tolerance.  In addition to computing these accuracy metrics, we record the overall integrator effort, as estimated by the number of time steps required at each of the fast and slow time scales.

\subsubsection{Controller robustness}
\label{sec:controller_robustness}

\begin{figure}[!th]
  \centering
  \includegraphics[trim={0 0 260 0}, clip, width=0.265\linewidth]{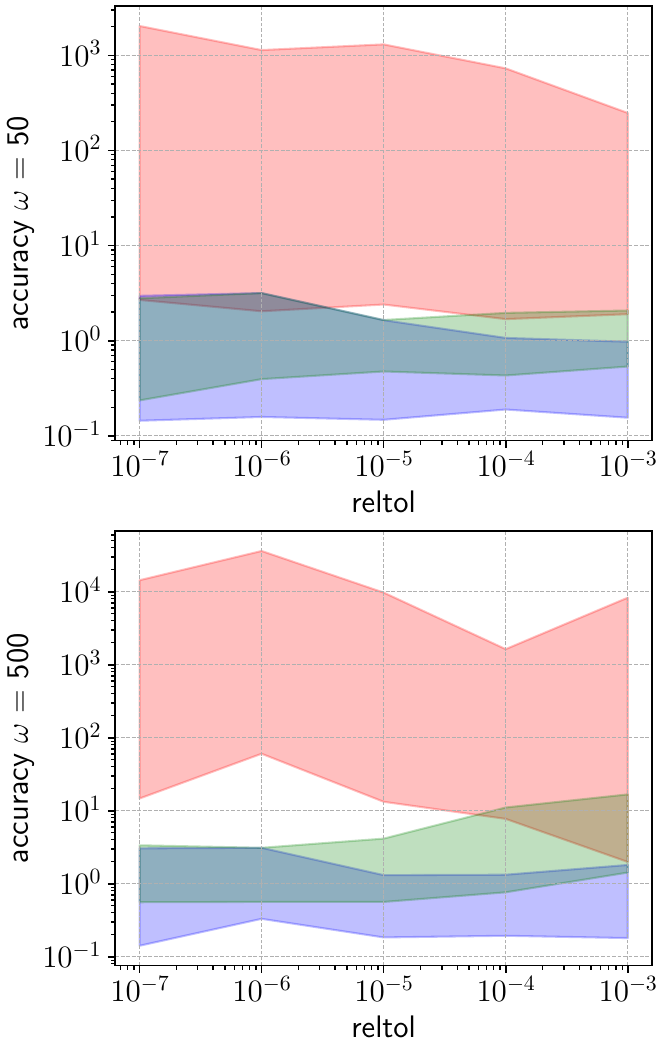}
  \includegraphics[trim={20 0 260 0}, clip, width=0.25\linewidth]{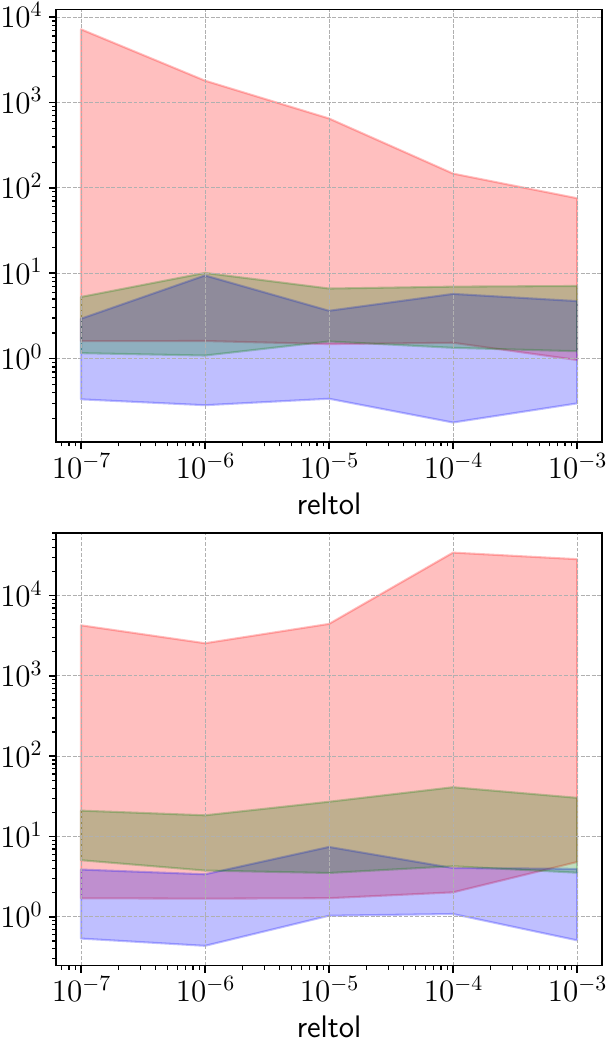}
  \includegraphics[trim={20 0 20 0}, clip, width=0.45\linewidth]{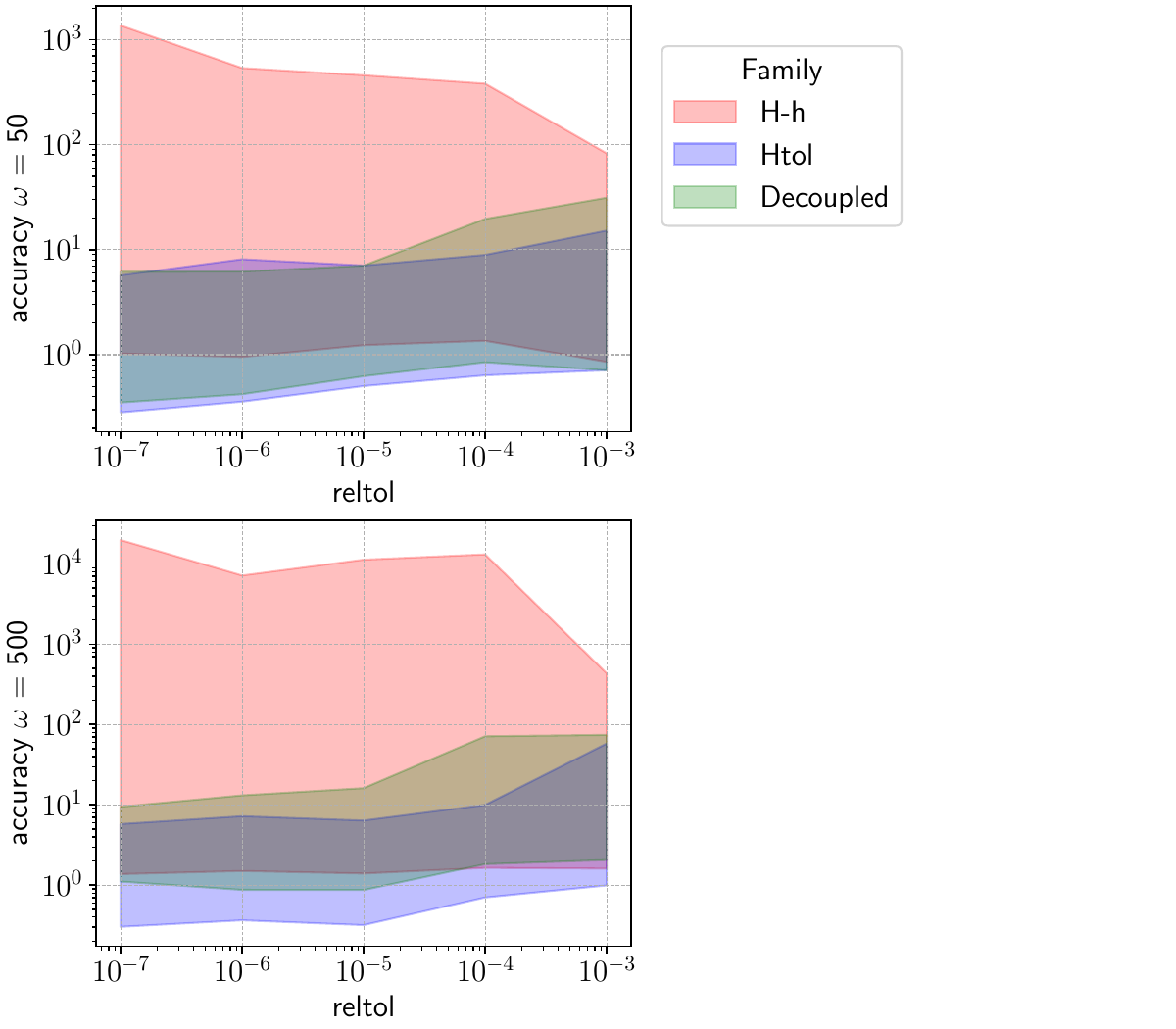}\\
  \includegraphics[trim={0 0 260 0}, clip, width=0.265\linewidth]{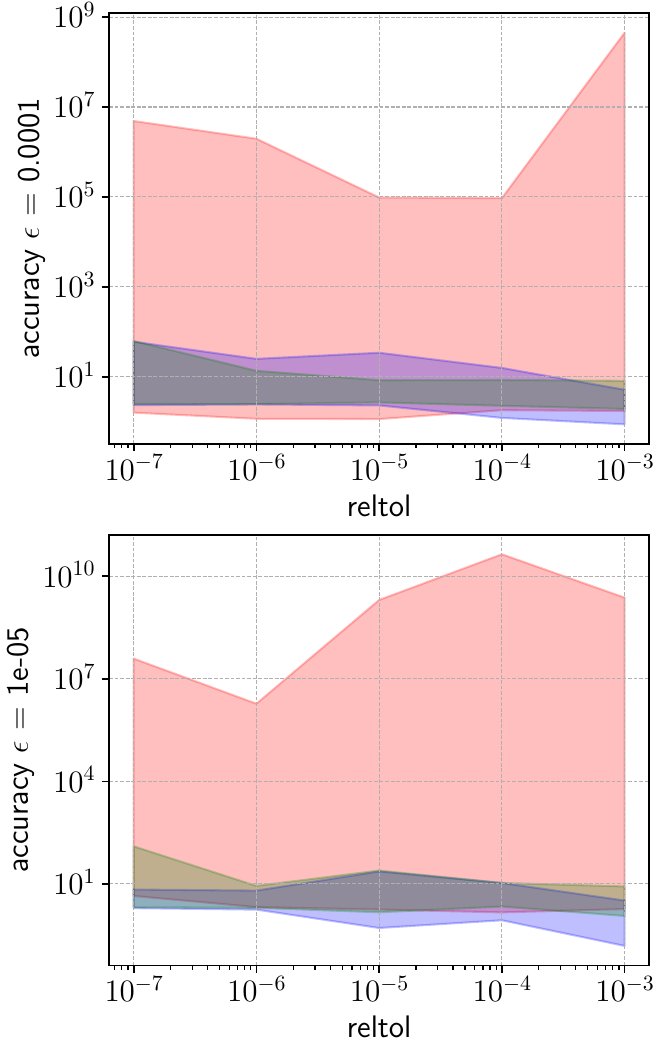}
  \includegraphics[trim={20 0 260 0}, clip, width=0.25\linewidth]{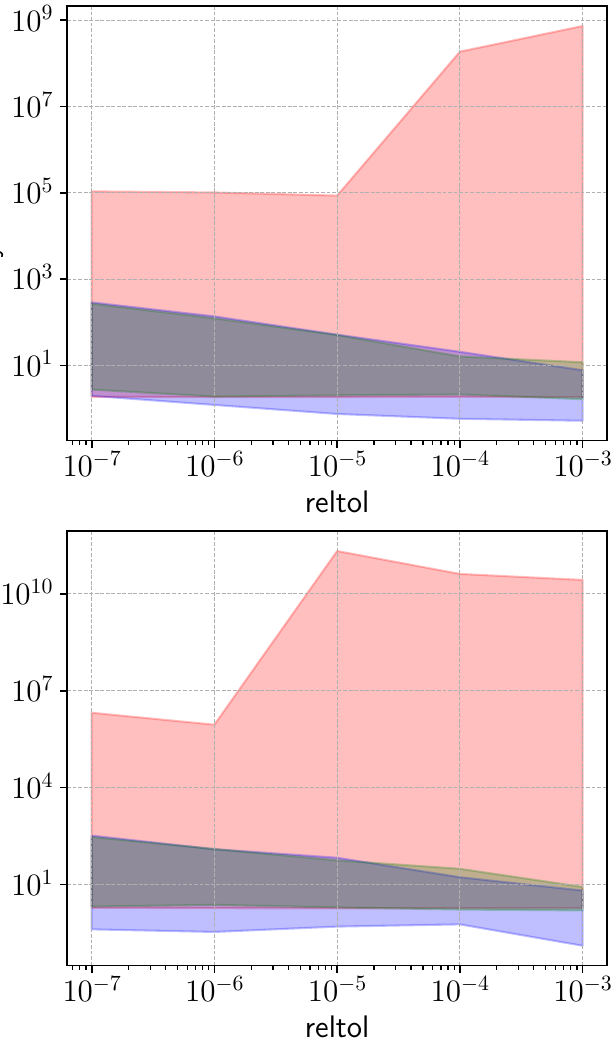}
  \includegraphics[trim={20 0 260 0}, clip, width=0.25\linewidth]{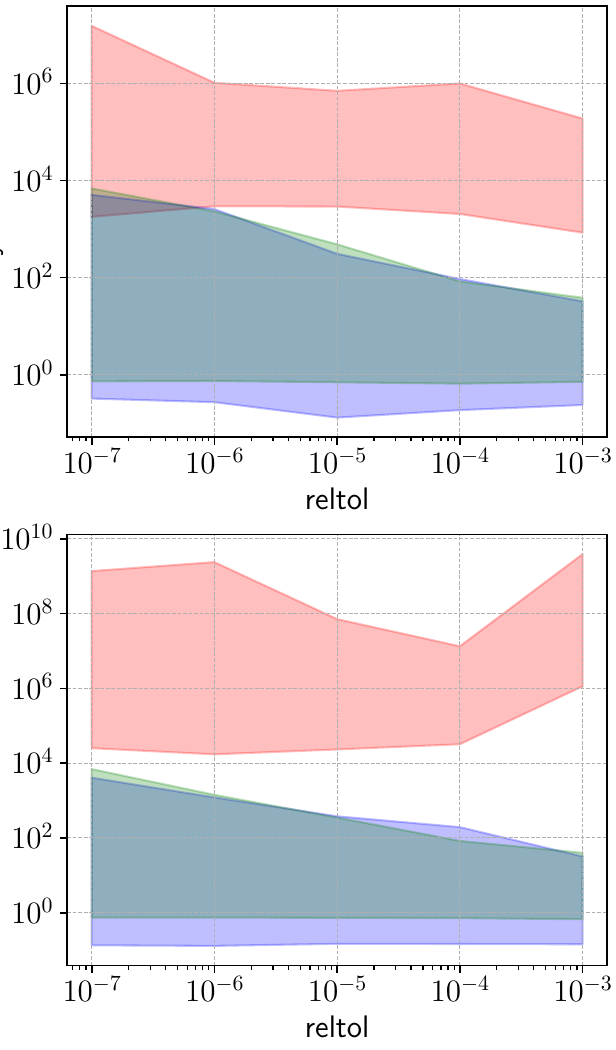}\hspace{2.6cm}

  \caption{Accuracy measurements for each multirate controller family, when tested across a wide range of tolerances and individual MRI controllers.  Columns denote MRI method accuracies: left are $\mathcal{O}(H^2)$, middle are $\mathcal{O}(H^3)$, and right are $\mathcal{O}(H^4)$ and $\mathcal{O}(H^5)$.  The KPR test problem is in the top two rows, and the Brusselator test problem is in the bottom two rows.  Note that where two regions overlap, their shading mixes together; the \HTol\ and \Decoupled\ families perfectly overlap, resulting in a teal-gray color.}
  \label{fig:shaded-accuracy}
\end{figure}

We begin with a simple question: how effective is each MRI adaptivity family (\Decoupled, \HTol, and \Hh) at achieving a range of target tolerances according to the metric \eqref{eq:accuracy}?  To answer this high-level question, Figure \ref{fig:shaded-accuracy} shows aggregated results across the full suite of tests.  Here, the top six plots show the KPR problem having multirate ratios $\omega=\{50,500\}$ for second order MRI methods (left), third order MRI methods (middle), and higher-order MRI methods (right), and the bottom six plots show similar results for the Brusselator problem with stiffness factors $\epsilon=\{10^{-4},10^{-5}\}$.  Within each plot, we present shaded regions for each family, showing the range of ``accuracy'' values attained by MRI methods and controllers within that family for each tolerance.  We can assess the robustness of each family by examining how tightly these plots cluster around the target accuracy of one.  Due to their disparate results, we separate our comments on the \Decoupled\ and \HTol\ controller families from the \Hh\ family.

For both problems and multirate regimes, the \Decoupled\ and \HTol\ controller families performed excellently, typically achieving solutions within a factor of 10 from the requested tolerances.  There were a few outliers among these for the KPR problem at $\omega=500$ and the loosest relative tolerances $10^{-4}$ and $10^{-3}$, where for some decoupled controllers the \texttt{ESDIRK46a} and \texttt{ERK45a} methods had errors over 100x larger than requested.  Due to its stiffness, the Brusselator problem had slightly more outliers, with \texttt{ESDIRK34a} giving excessive error for some tighter tolerances, and \texttt{ESDIRK46a} struggling to attain solutions within 100x of the target across a wide range of tolerances and single-rate controllers.  We additionally note that from the plots in Figure \ref{fig:shaded-accuracy} it is difficult to discern significant performance differences between the \Decoupled\ and \HTol\ families, indicating that their accuracy depends more strongly on the underlying adaptive MRI method or single-rate controller than the time step selection mechanism.  Based on these results, we conclude that both of these families are robust across a wide range of MRI methods and tolerances, and should be applicable to most multirate applications.

The \Hh\ controller family did not fare as well.  Although some methods and controllers were able to meet the requested tolerances, most had errors that were orders of magnitude larger than requested.  Generally, the \Hh\ controllers deteriorated as the problems became more multirate; in separate tests on weakly multirate problems (e.g., KPR with $\omega=5$ and the stiff Brusselator with $\epsilon=0.01$) this family was generally able to meet the requested error tolerances, which agrees with the results from \cite{fishAdaptiveTimeStep2023}.  The only outliers from this generally poor performance were the \texttt{MERK21} and \texttt{MERK32} methods, that achieved the target accuracy for all tolerances and multirate parameters.  \emph{No other combination of adaptive MRI method and \Hh\ controller was able to achieve even within 10x of the target accuracy for the challenging problems with $\omega=500$ and $\epsilon=10^{-5}$}.  We cannot yet explain why the embedded MERK methods with \Hh\ controllers outperform other embedded MRI methods; however, \corrA{these results reinforce our analysis from Section \ref{subsubsec:Hh-analysis} showing that the \Hh\ controllers struggle in comparison with the \Decoupled\ and \HTol\ families for applications with significant multirate nature.}

\subsubsection{Adaptive MRI efficiency}
\label{sec:adaptive_efficiency}

We now turn our attention to the computational efficiency of the individual adaptive MRI methods.  To test computational efficiency for a given embedded MRI method and multirate controller, we ran each problem using relative solution tolerances $10^{-k}, k=3,4,\ldots,7$.  As in \cite{chinomonaImplicitExplicitMultirateInfinitesimal2021,fishAdaptiveTimeStep2023,luanNewClassHighOrder2020,luanMultirateExponentialRosenbrock2022}, we then estimate the computational ``work'' required for a calculation by separately considering the number of slow steps and the total number of fast and slow steps.  The former is relevant for multirate applications where the slow operators require significantly more computational effort than the fast operators, and thus MRI methods are applied to reduce the number of slow time steps.  The latter is relevant for applications where the slow and fast operators require similar computational effort.  For either metric, we then plot the computational efficiency by graphing the solution error as a function of work, overlaying these curves for a variety of methods.  For any desired accuracy, the left-most curve thus corresponds with the most efficient method.

In lieu of presenting plots that overlay hundreds of combinations of MRI methods and adaptive controllers, we first downselected only the most performant combinations.  Within each test problem (KPR with $\omega=\{50,500\}$, and stiff Brusselator with $\epsilon=\{10^{-4},10^{-5}\}$), MRI method order (2, 3, and higher), and work metric (slow steps, fast steps), we:
\begin{itemize}
\item define a set of 20 logarithmically-spaced test tolerances in $[10^{-6},10^{-2}]$;
\item at each test tolerance, we estimate the work required to attain the target tolerance by interpolating the (work, error) pairs from our data;
\item rank each method+controller combination for each test tolerance, accumulating their average rank over all test tolerances.
\end{itemize}

To more rigorously compare the performance of MRI methods and controller combinations, we also conducted a variety of statistical analyses of the full set of efficiency data.  We used a one-way Analysis of variance (ANOVA) to analyze controllers and a repeated measure ANOVA for each MRI method, examining fast and slow time scales separately. Both analyses yielded p-values of zero, indicating a statistically significant difference between controllers and between methods of a given order.

To identify the best-performing method of a given order for either the fast or slow work metric on each test problem, we first grouped the data by controller. For each controller and test problem, we calculated the mean and standard deviation of all rank values across methods of that particular order. Next, we condensed each method's performance into a single value by averaging its rank values.  Finally, we calculated the z-score for each method to determine its deviation from the mean, and then averaged these z-scores across all controllers. An MRI method with negative z-score (i.e., below the mean) indicates better performance, while a positive z-score (above the mean) signifies poorer performance. A similar z-score analysis is conducted to determine the best-performing controllers across all methods and test problems, for the fast and slow time scales.  We include these statistical results in our discussion below.

\textbf{Family performance.}
As observed in our robustness tests in Section \ref{sec:controller_robustness}, the largest determining factor for adaptive efficiency was the controller family; this was followed by the embedded MRI method itself, and then the specific choice of controller.  Specifically, when comparing the average rankings of each family across all test problems, MRI methods and work metrics, the \HTol\ and \Decoupled\ families had z-scores of -0.349 and -0.340, respectively, while the \Hh\ family had a z-score of 0.862.  This indicates that the proposed families performed similarly to each other, but were far more efficient than the \Hh\ family.

\begin{figure}[!tbh]
\centering
\subfloat[]{\label{}\includegraphics[width=.49\linewidth]{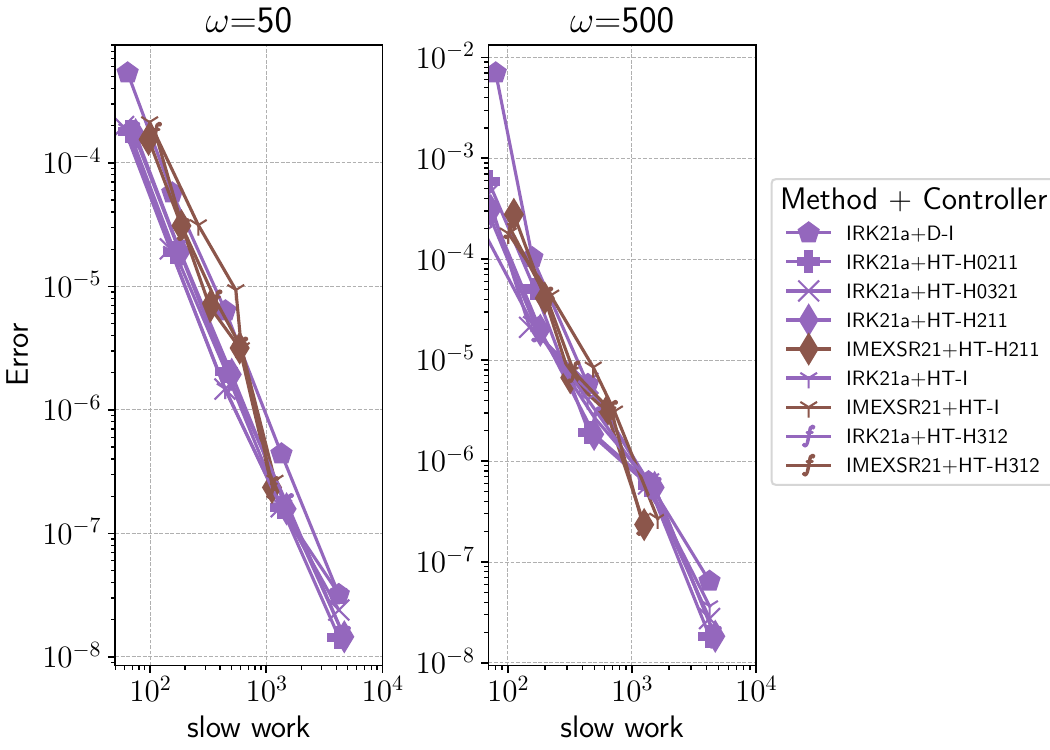}}\hfill
\subfloat[]{\label{}\includegraphics[width=.49\linewidth]{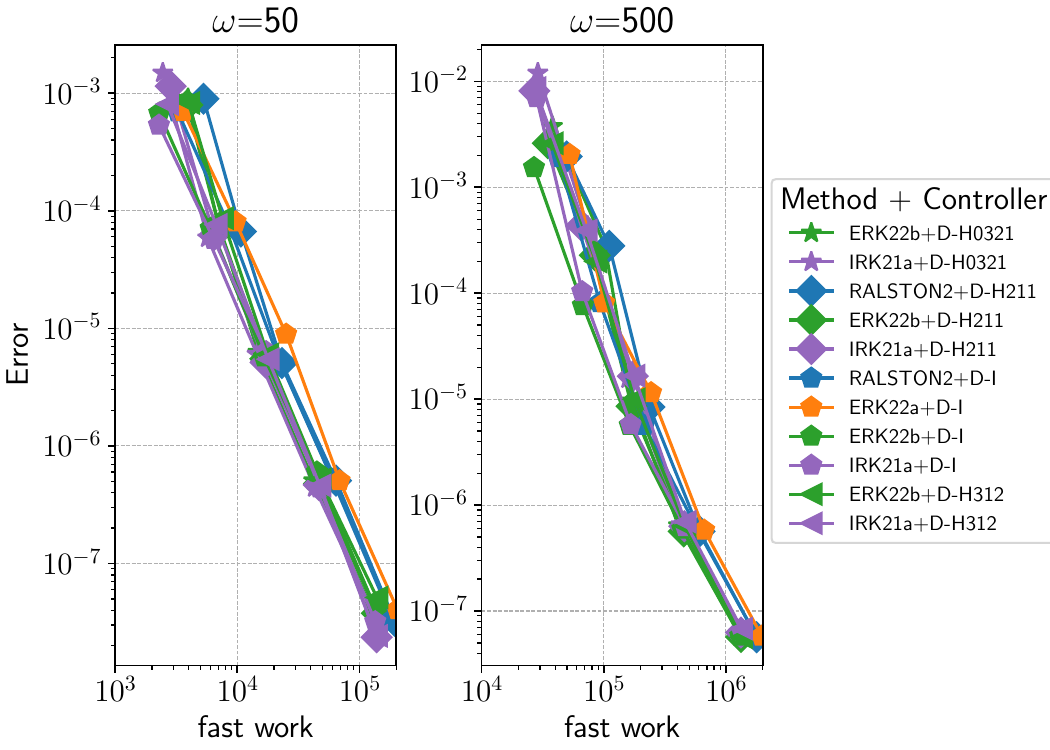}}\par
\subfloat[]{\label{}\includegraphics[width=.49\linewidth]{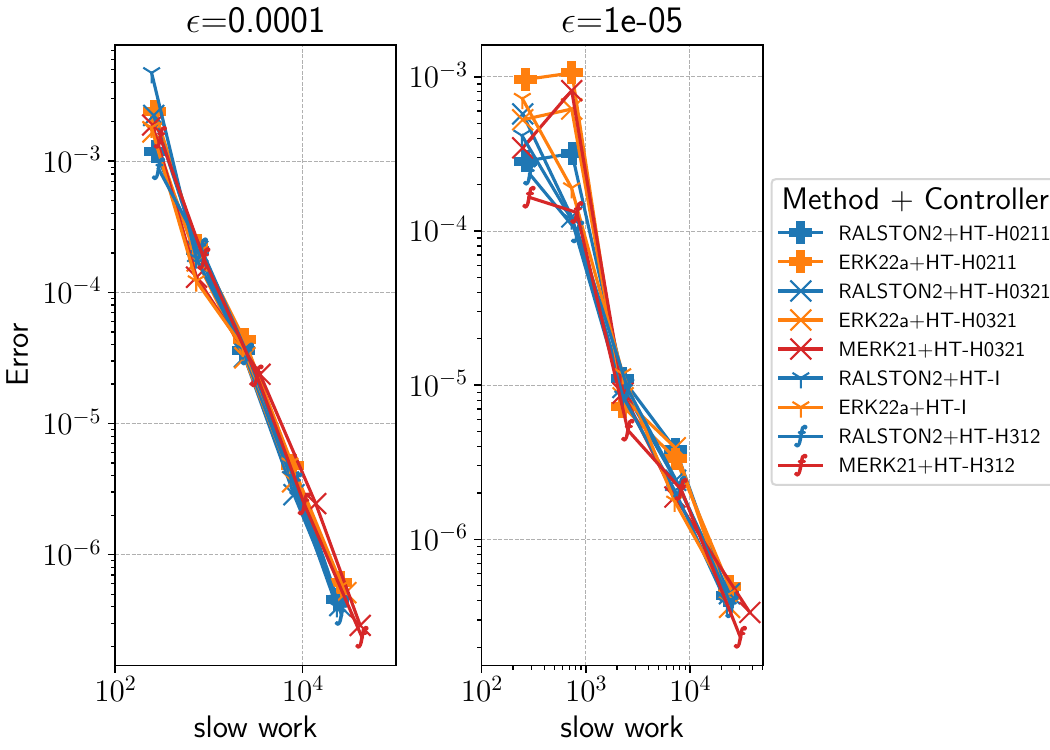}}\hfill
\subfloat[]{\label{}\includegraphics[width=.49\linewidth]{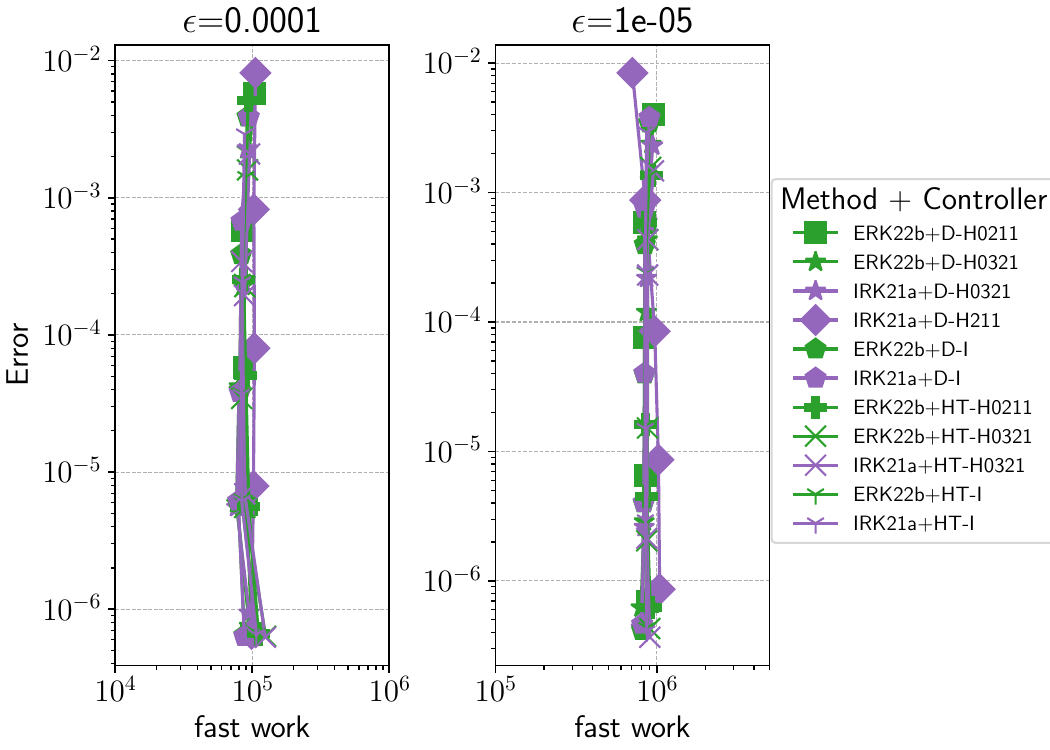}}\par

\caption{Efficiency comparisons for the top second-order adaptive MRI methods.  The top row contains the slow and fast time scales for the KPR test problem with multirate ratios $\omega=\{50,500\}$. The stiff Brusselator test problem with both stiffness parameters $\epsilon = \{10^{-4},10^{-5}\}$ is on the bottom.}
\label{fig:efficiency-lo}
\end{figure}

\begin{figure}[htbp!]
\centering
\subfloat[]{\label{}\includegraphics[width=.49\linewidth]{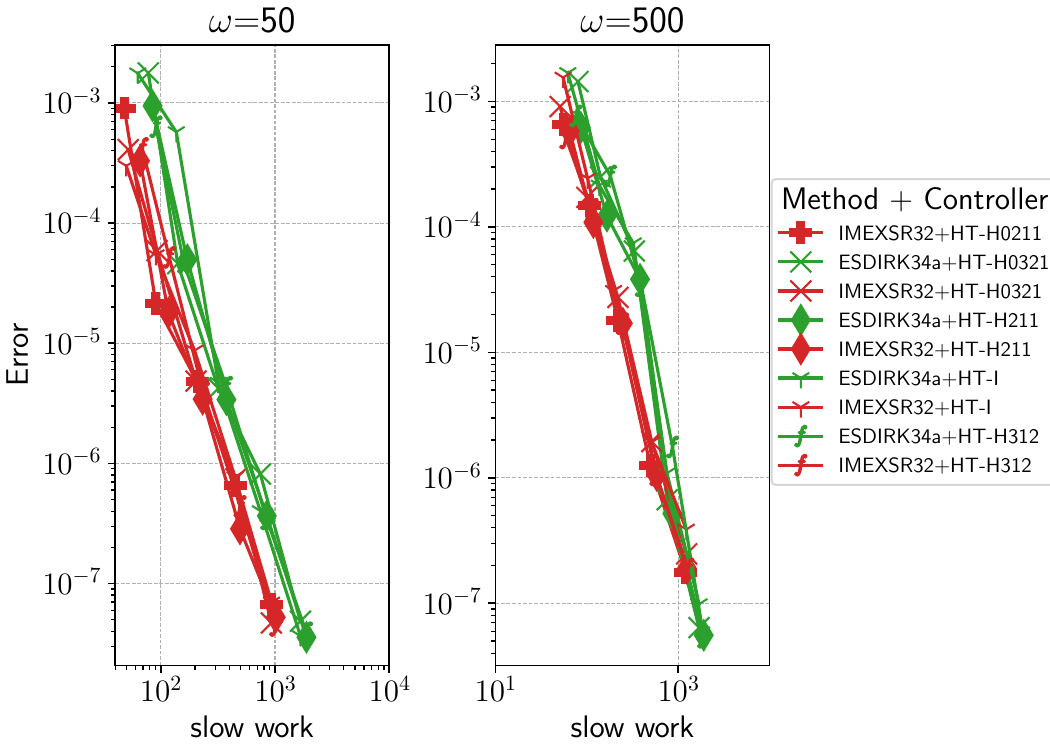}}\hfill
\subfloat[]{\label{}\includegraphics[width=.49\linewidth]{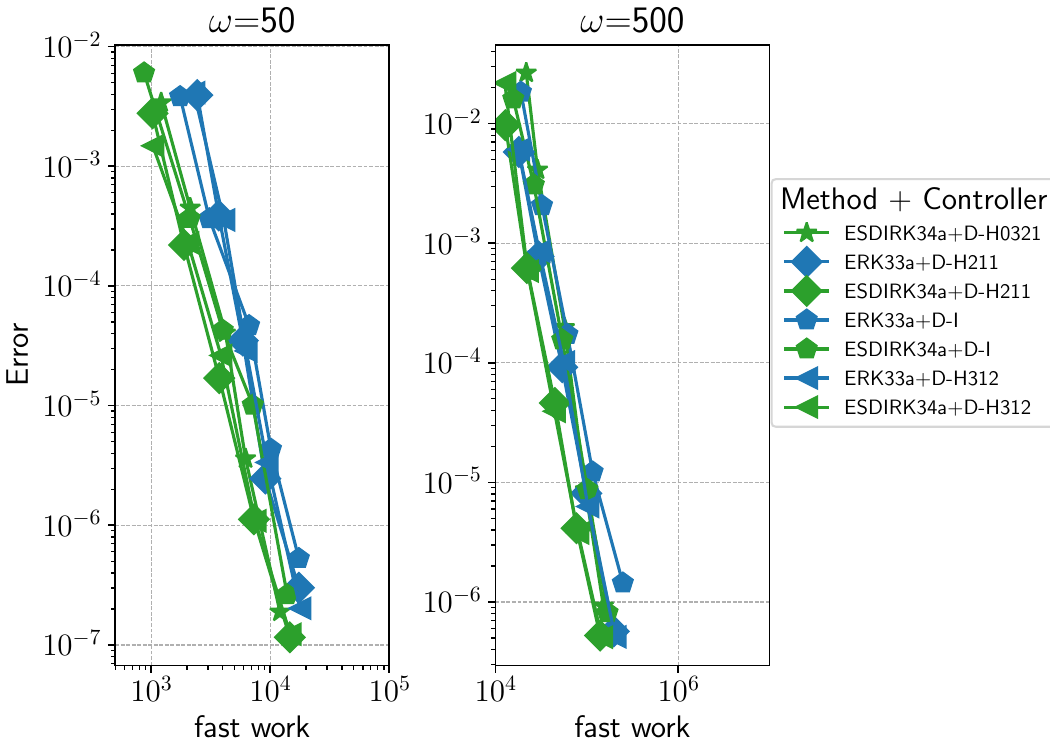}}\par
\subfloat[]{\label{}\includegraphics[width=.49\linewidth]{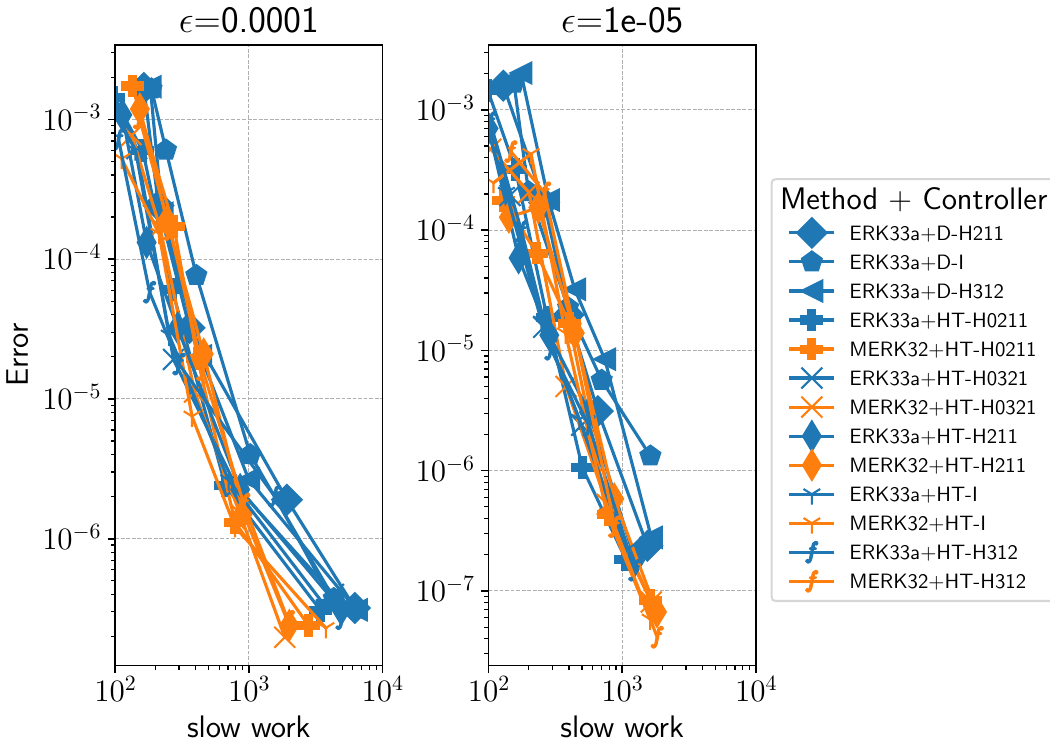}}\hfill
\subfloat[]{\label{}\includegraphics[width=.49\linewidth]{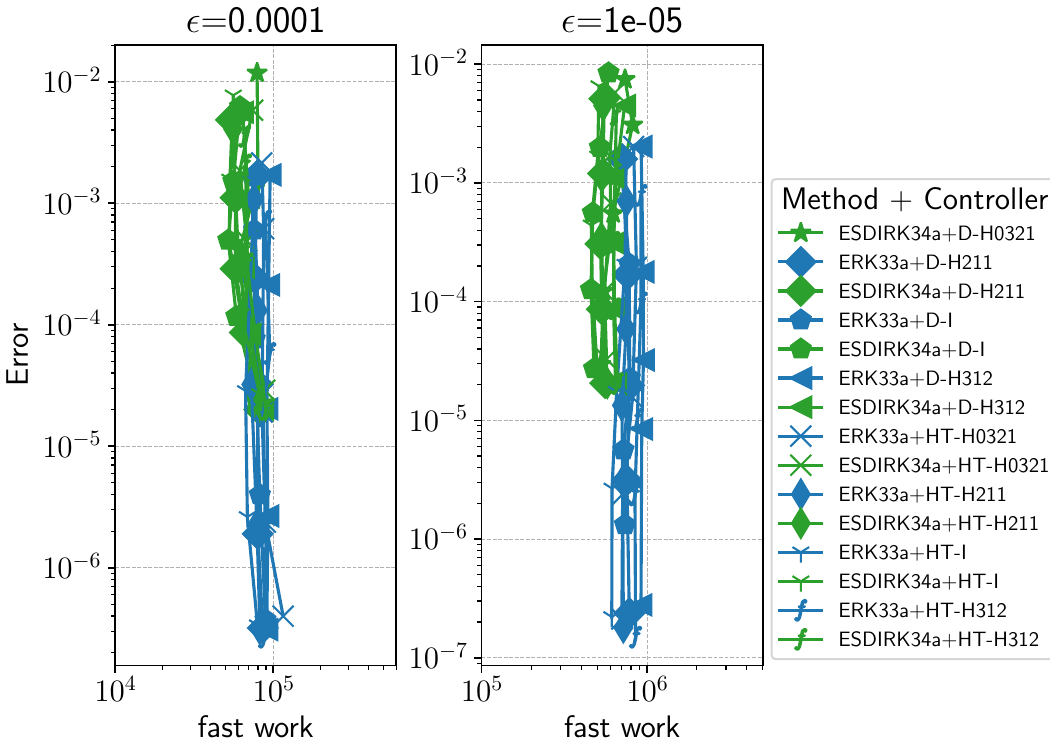}}\par

\caption{Efficiency comparisons for the top third-order adaptive MRI methods. The top row shows the slow and fast time scales for the KPR test problem with multirate ratios $\omega=\{50,500\}$. The stiff Brusselator test problem with both stiffness parameters $\epsilon = \{10^{-4},10^{-5}\}$ is on the bottom row.}
  \label{fig:efficiency-mid}
\end{figure}

\begin{figure}[!tbh]
\centering
\subfloat[]{\label{}\includegraphics[width=.49\linewidth]{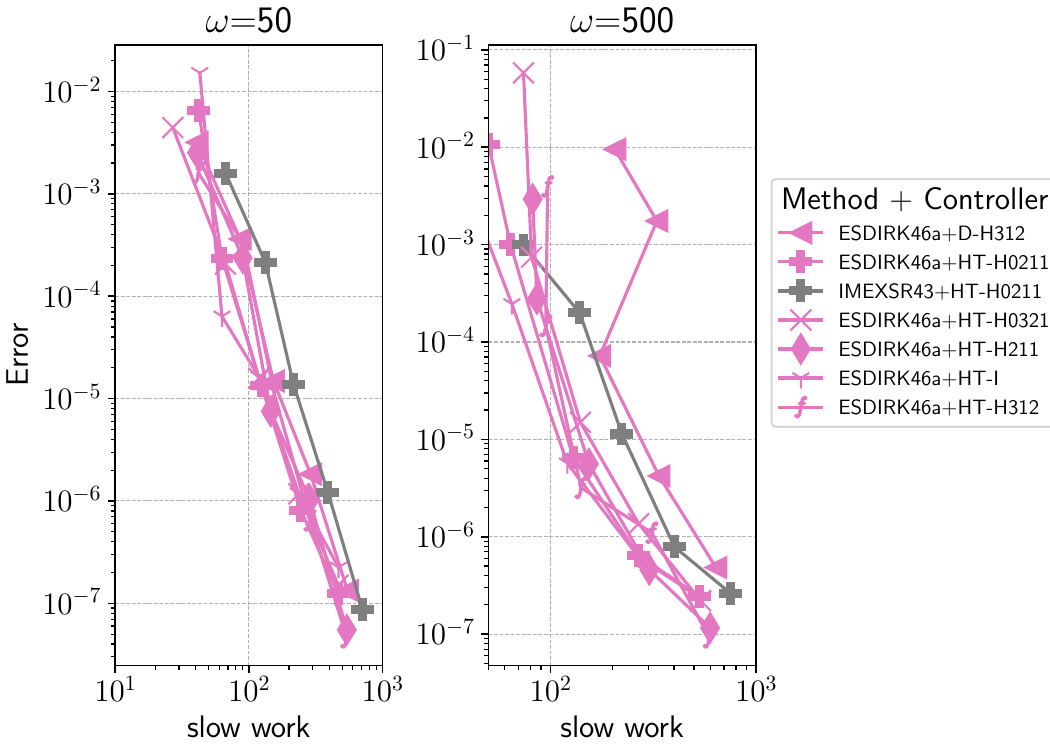}}\hfill
\subfloat[]{\label{}\includegraphics[width=.49\linewidth]{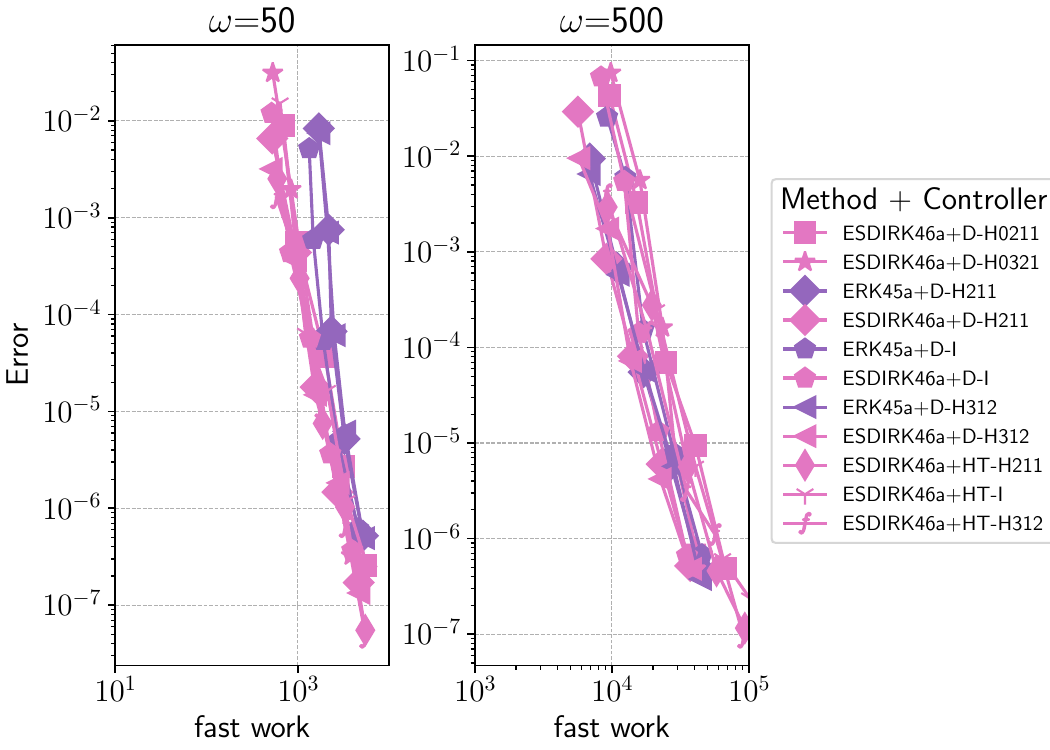}}\par
\subfloat[]{\label{}\includegraphics[width=.49\linewidth]{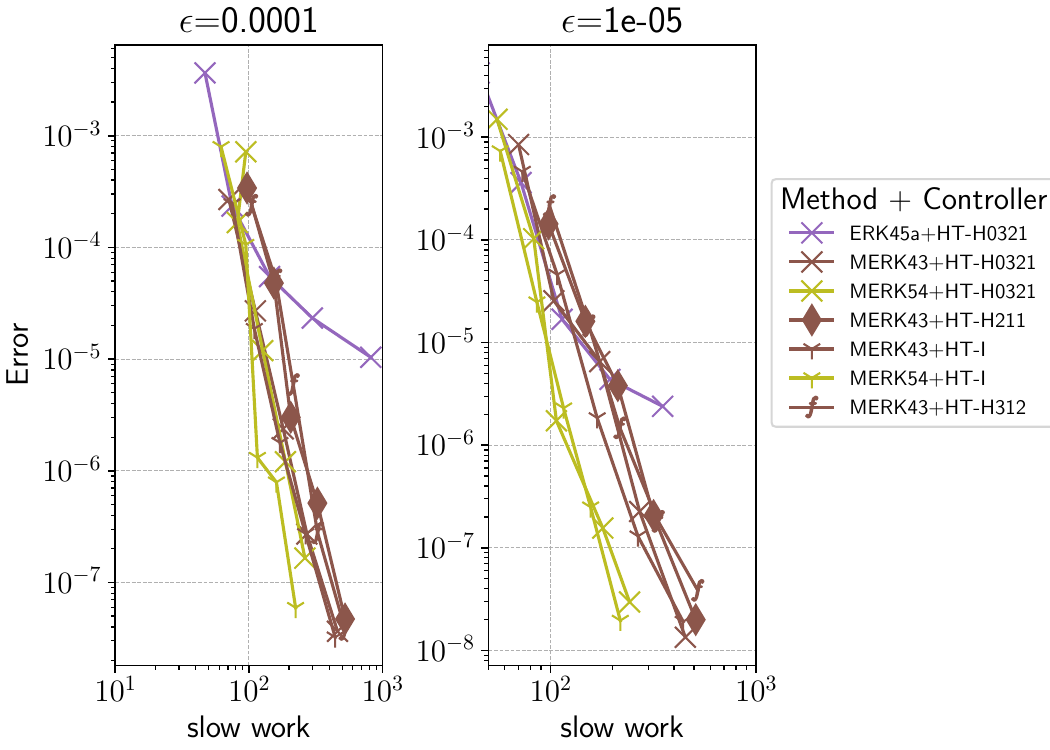}}\hfill
\subfloat[]{\label{}\includegraphics[width=.49\linewidth]{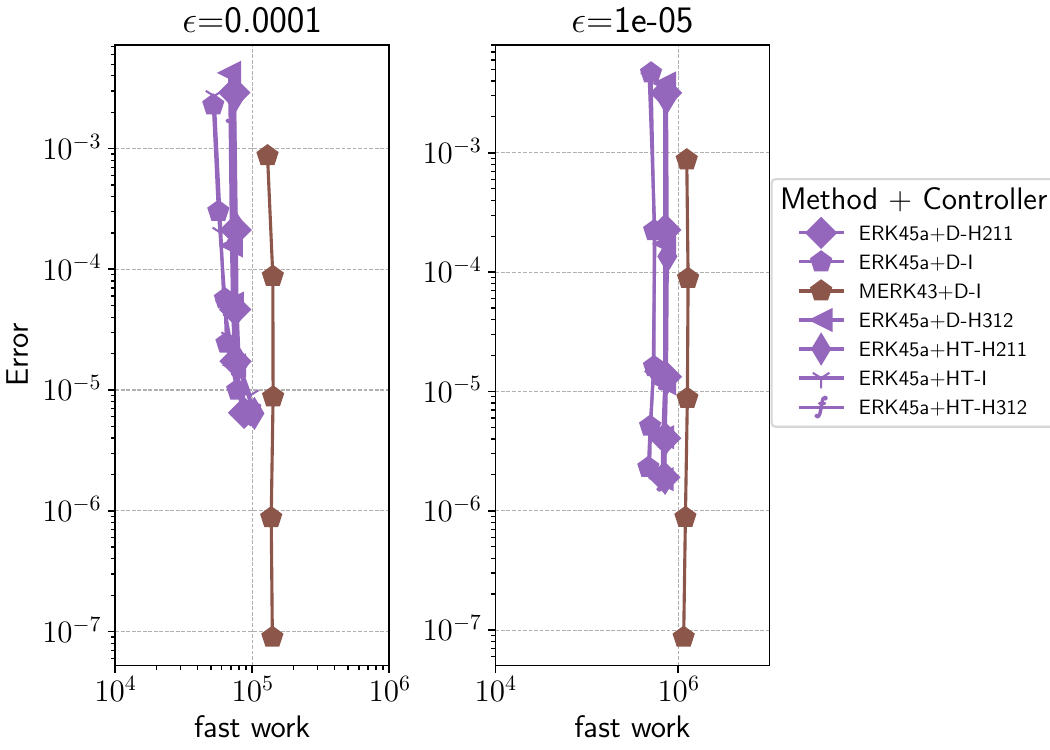}}\par

\caption{Efficiency comparisons for the top fourth- and fifth-order adaptive MRI methods.  The top row shows the slow and fast time scales for the KPR test problem with multirate ratios $\omega=\{50,500\}$. The stiff Brusselator test problem with both stiffness parameters $\epsilon = \{10^{-4},10^{-5}\}$ is on the bottom row.}
  \label{fig:efficiency-hi}
\end{figure}

\textbf{MRI method performance.}
To compare individual embedded MRI methods, we downselected to the twelve fastest MRI method and controller pairs for each combination of problem, method order, and work metric. This process yielded eight sets of top-performing pairs.  We then computed the intersection of the two sets corresponding to the different multirate values, retaining only those pairs that consistently rank among the top 12 for both. In Figures \ref{fig:efficiency-lo}-\ref{fig:efficiency-hi}, we overlay the efficiency plots for these sets of ``fastest'' methods to determine the most performant MRI method and adaptivity controller pairs.  In each figure, the KPR test problem is on the top row, and the stiff Brusselator test problem is on the bottom rows.  Within each row, we show the plots for slow computational efficiency on the left.

\begin{table}
  \centering
  \caption{Average rank z-scores for embedded second-order MRI methods.}
  \begin{tabular}{l|rr|r|rr|r|}
    \small
    MRI    & \multicolumn{3}{|c|}{Slow} & \multicolumn{3}{|c|}{Fast} \\
    Method & KPR & Bruss & Avg & KPR & Bruss & Avg \\
    \hline
    \texttt{IRK21a}   & -1.50 & -0.17 & -0.84 & -0.95 & -1.09 & -1.02 \\
    \texttt{RALSTON2} &  0.48 & -0.76 & -0.14 & -0.35 & -0.17 & -0.26 \\
    \texttt{ERK22a}   &  0.32 & -0.34 & -0.01 &  0.10 &  0.16 &  0.13 \\
    \texttt{MERK21}   &  0.43 & -0.09 &  0.17 &  0.82 &  0.76 &  0.79 \\
    \texttt{ERK22b}   &  0.86 & -0.47 &  0.19 & -1.11 & -1.15 & -1.13 \\
    \texttt{IMEXSR21} & -0.58 &  1.84 &  0.63 &  1.49 &  1.49 &  1.49 \\
    \hline
  \end{tabular}
  \label{tab:z-scores-lo}
\end{table}

Figure \ref{fig:efficiency-lo} shows the efficiency results for the best second-order adaptive MRI method combinations.  The top MRI methods showed the most variability in their performance at the slow time scale (left plots).  Here, different MRI methods excelled for each problem, with \texttt{IRK21a} the best for the KPR problem across a range of \HTol\ controllers (plus one \Decoupled\ controller), while \texttt{RALSTON2} and \texttt{ERK22a} dominated the stiff Brusselator problem (again using various \HTol\ controllers).  The corresponding z-scores for each MRI method's average rank are given in Table \ref{tab:z-scores-lo}.  At the fast time scale (right plots), the MRI methods showed more consistent performance across both test problems, with \texttt{ERK22b} and \texttt{IRK21a} giving the best performance, and \texttt{IMEXSR21} the worst performance.  Interestingly, when measuring fast cost the \Decoupled\ controllers rank near the top more frequently than \HTol.  We additionally note a feature that will be present in higher-order methods as well -- the fast time scale stiff Brusselator work-precision curves are nearly vertical, a characteristic of explicit methods on stiffness-dominated problems.  The z-scores for these methods are also given in Table \ref{tab:z-scores-lo}.  From these results, we conclude that the second-order adaptive MRI methods are generally effective at both fast and slow time scales, with \texttt{IRK21a} and \texttt{RALSTON2} being the most efficient across both problems and time scales, while \texttt{IMEXSR21} was generally the least efficient second-order method.

\begin{table}
  \centering
  \caption{Average rank z-scores for embedded third-order MRI methods.}
  \begin{tabular}{l|rr|r|rr|r|}
    \small
    MRI    & \multicolumn{3}{|c|}{Slow} & \multicolumn{3}{|c|}{Fast} \\
    Method & KPR & Bruss & Avg & KPR & Bruss & Avg \\
    \hline
    \texttt{ERK33a}    &  0.10 & -1.01 & -0.46 & -0.30 & -0.86 & -0.58\\
    \texttt{ESDIRK34a} &  0.05 &  0.21 &  0.13 & -1.25 & -0.83 & -1.04\\
    \texttt{IMEXSR32}  & -1.10 &  1.38 &  0.14 &  1.17 &  1.36 &  1.27\\
    \texttt{MERK32}    &  0.95 & -0.58 &  0.19 &  0.37 &  0.33 &  0.35\\
    \hline
  \end{tabular}
  \label{tab:z-scores-mid}
\end{table}

Figure \ref{fig:efficiency-mid} shows the efficiency of the best third-order adaptive MRI method combinations.  Again, for problems where the cost is dominated by operators at the slow time scale, the KPR and stiff Brusselator problems have different optimal MRI methods.  For KPR, by far the most efficient MRI method was \texttt{IMEXSR32} for an array of \HTol\ controllers, while \texttt{ERK33a} and \texttt{MERK32} are the most efficient for the stiff Brusselator (again predominantly using \HTol\ controllers).  Also similarly to the second-order methods, for multirate applications where the slow and fast operators have commensurate cost, both the KPR and stiff Brusselator tests agree that the most efficient third-order methods are \texttt{ERK33a} and \texttt{ESDIRK34a}, with the \Decoupled\ controllers appearing more frequently than \HTol.  Interestingly, we note that for the stiff Brusselator problem at the fast time scale, although \texttt{ESDIRK34a} is more efficient at loose tolerances, it is incapable of achieving accuracies better than $10^{-5}$.  Table \ref{tab:z-scores-mid} presents the z-scores for each problem and metric, where it is clear that the overall top-performing third-order MRI method was \texttt{ERK33a}, since its average z-scores for both slow and fast work metrics were negative.

\begin{table}
  \centering
  \caption{Average rank z-scores for embedded fourth- and fifth-order MRI methods.}
  \begin{tabular}{l|rr|r|rr|r|}
    \small
    MRI    & \multicolumn{3}{|c|}{Slow} & \multicolumn{3}{|c|}{Fast} \\
    Method & KPR & Bruss & Avg & KPR & Bruss & Avg \\
    \hline
    \texttt{ERK45a}    &  0.16 & -0.84 & -0.34 & -0.59 & -1.28 & -0.93\\
    \texttt{ESDIRK46a} & -1.34 &  1.11 & -0.11 & -1.20 &  0.23 & -0.48\\
    \texttt{MERK43}    &  0.95 & -0.89 &  0.03 &  0.49 & -0.15 &  0.17\\
    \texttt{MERK54}    &  0.44 & -0.20 &  0.12 & -0.13 & -0.26 & -0.20\\
    \texttt{IMEXSR43}  & -0.21 &  0.82 &  0.30 &  1.43 &  1.46 &  1.44\\
    \hline
  \end{tabular}
  \label{tab:z-scores-hi}
\end{table}

We present efficiency results for fourth- and fifth-order adaptive MRI methods in Figure \ref{fig:efficiency-hi} and Table \ref{tab:z-scores-hi}.  Here, we see that the best-performing methods again differ between the KPR and stiff Brusselator problems.  For the KPR problem at both $\omega$ values, the implicit \texttt{ESDIRK46a} provides the best efficiency using either work metric, while \texttt{ERK45a} is competitive when considering fast time scale effort.  For the stiff Brusselator problem, we see that only \texttt{MERK43} and \texttt{MERK54} are able to achieve solutions with errors below $10^{-6}$, but that for looser tolerances \texttt{ERK45a} is competitive at the slow time scale and wins at the fast time scale.  We additionally note that as with the lower order methods, when computational cost is dominated by the slow time scale operators, nearly all of the top-performing methods use \HTol\-based controllers, while for problems where fast scale computational effort is significant the \Decoupled\ controllers slightly outnumber \HTol.

\textbf{Multirate controller performance.}
We conclude this section by turning our attention to the performance of individual multirate controllers.  Our previous statistical analyses compared MRI methods when using the same multirate controllers, so here we compare the performance of each \HTol\ and \Decoupled\ controller when using the same MRI methods.  The z-scores from this analysis are shown in Table \ref{tab:z-scores-controllers}.  Here, we perform the analysis separately for the average ranks at the slow time scale and the fast time scale in the left and right columns, respectively.  These underscore the intuitive results seen above, that the \HTol\ controllers are clearly more efficient than the \Decoupled\ controllers when slow time scale costs are dominant, but that \Decoupled\ are more efficient when the fast time scale costs become significant.  These further indicate that within each family, the multirate controllers constructed from the single-rate ``I'' controller are significantly more efficient for these test problems than the others.

\begin{table}[!th]
  \centering
  \caption{Average rank z-scores for multirate controllers (compares only the \HTol\ and \Decoupled\ families).  Ranked scores for slow time scale work are on the left, and
  for fast time scale work are on the right.}
  \begin{tabular}{l r || l r}
    \multicolumn{2}{c||}{Slow Scale} & \multicolumn{2}{c}{Fast Scale}\\
    Multirate controller & z-score & Multirate controller & z-score\\
    \hline
    HT-I      &  -0.690  &  D-I       &   -0.557 \\
    HT-H0321  &  -0.349  &  D-H312    &   -0.419 \\
    HT-H312   &  -0.280  &  D-H211    &   -0.405 \\
    HT-H0211  &  -0.264  &  D-H0321   &   -0.181 \\
    HT-H211   &  -0.255  &  D-H0211   &    0.037 \\
    D-I       &   0.065  &  HT-I      &    0.150 \\
    D-H312    &   0.307  &  HT-H211   &    0.189 \\
    D-H211    &   0.420  &  HT-H312   &    0.245 \\
    D-H0321   &   0.445  &  HT-H0321  &    0.443 \\
    D-H0211   &   0.601  &  HT-H0211  &    0.499 \\
    \hline
  \end{tabular}
  \label{tab:z-scores-controllers}
\end{table}

\subsubsection{Temporal Adaptivity Comparisons}
\label{sec:adaptivity_comparisions}

\begin{figure}[!th]
    \centering
    \includegraphics[width=0.9\linewidth]{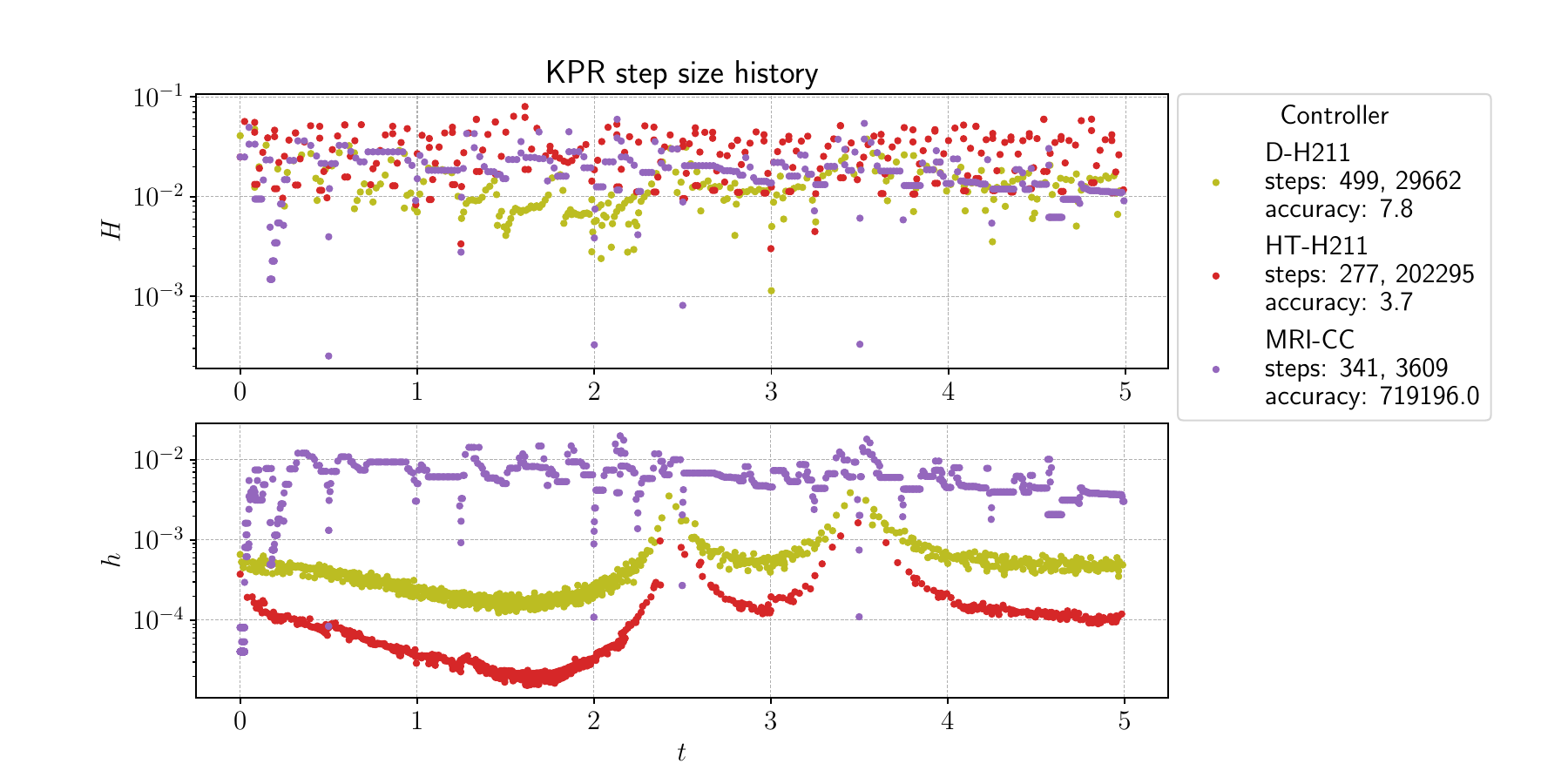}\\
    \includegraphics[width=0.9\linewidth]{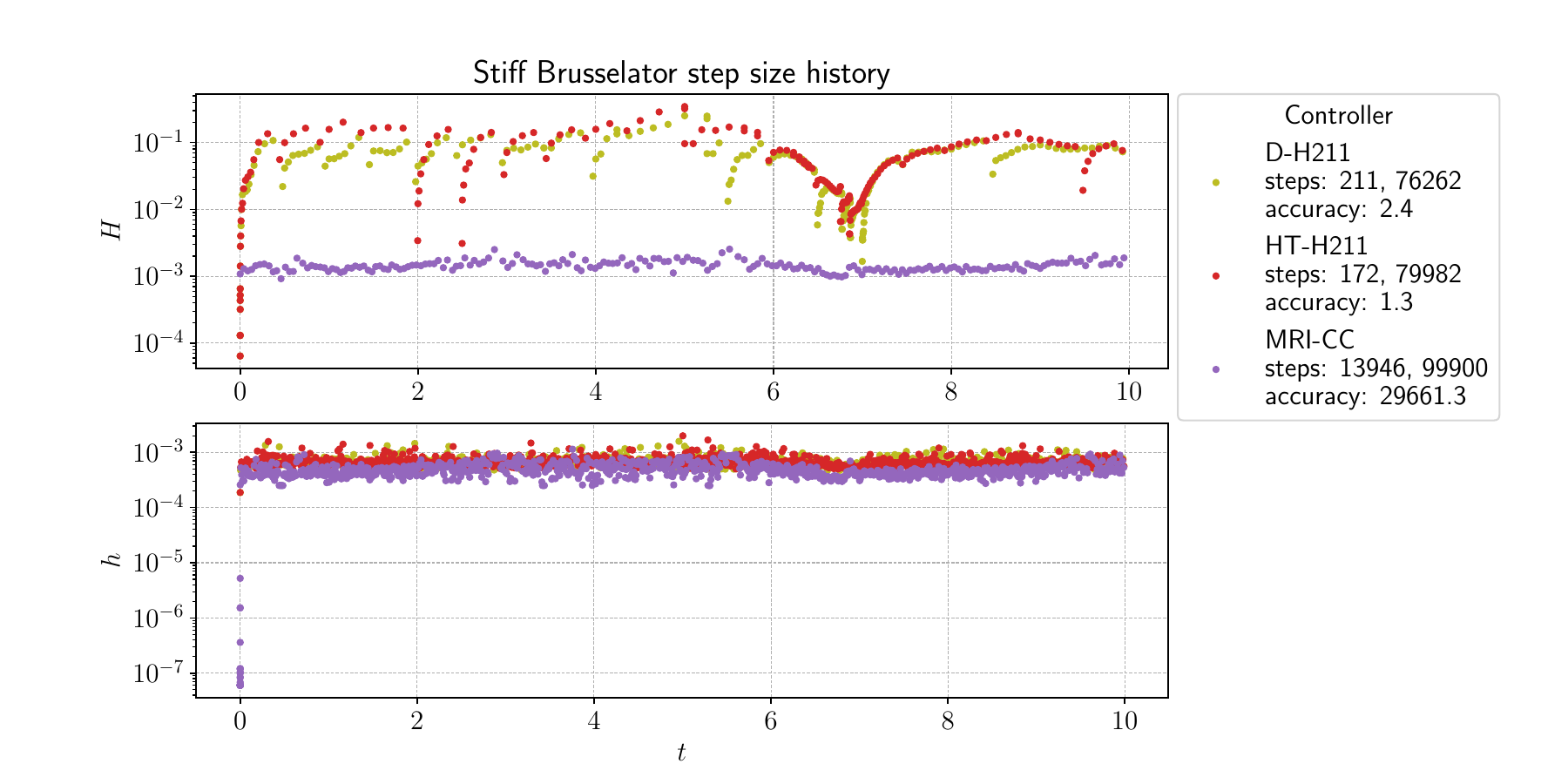}

    \caption{Slow and fast time step histories for each multirate controller family on the KPR and stiff Brusselator test problems.  The legends list the numbers of slow and fast time steps, as well as the attained solution accuracy factor \eqref{eq:accuracy}.}
    \label{fig:adaptivity-comparisons}
\end{figure}

To better understand the differing behavior between each controller family, we present time histories of $H$ and $h$ as a function of the internal simulation time, $t$, on the KPR and stiff Brusselator tests.  To focus the presentation, all simulations use absolute and relative tolerances $10^{-11}$ and $10^{-4}$, the \texttt{ERK33a} MRI method, and ARKODE's default third-order explicit Runge--Kutta method at the fast time scale.  The KPR problem is run with $e_s=5$, $e_f=0.5$, and $\omega = 500$, and the stiff Brusselator problem uses $\epsilon = 10^{-4}$.  With this setup, we selected the D-H211, HT-H211, and MRICC multirate controllers.  To help unclutter the figures, we plot a subset of these internal step sizes in Figure \ref{fig:adaptivity-comparisons}: up to 200 values of $H(t)$ and up to 1000 values of $h(t)$.  In the legend, we list the numbers of slow and fast time steps, and the attained ``accuracy'' as defined in equation \eqref{eq:accuracy}.  As expected from the analytical solution to the KPR problem, we see that at $t \approx 2.5$ and $t\approx 3.5$ the oscillations in the fast variable slow down, allowing the fast integrator to increase its steps dramatically.  Similarly, due to rapid changes in the slow variables for the stiff Brusselator solution at $t\approx 6.5$ the ``slow'' solution components speed up rapidly, causing the slow integrator to shrink steps while the fast integrator remains steady.

From these plots, we additionally discern the underlying differences between each controller family.  Notably, while all methods exhibit similar slow step histories on the KPR problem, at the fast time scale the MRI-CC \Hh\ controller consistently took overly large fast time steps, resulting in solutions with orders of magnitude more error than requested.  Additionally, the KPR plots show the tradeoff between the \Decoupled\ and \HTol\ controllers, namely that the \HTol\ controller will shift additional effort onto the fast time scale so that it can increase the slow step size.  The stiff Brusselator problem tells a slightly different story: although all controllers exhibited similar behavior when resolving the fast/stiff time scale, the MRI-CC controller took slow time steps that were far smaller than necessary, resulting in significantly worse efficiency (and interestingly, also much worse accuracy).  Meanwhile, the \HTol\ controller required marginally fewer slow steps than its \Decoupled\ counterpart.

\subsection{Nested Multirate Tests}
\label{sec:nested_tests}

We end with a demonstration of the multirate \HTol\ controller on a nested multirate application.  For this, we consider a modification of the previous KPR problem \eqref{eq:KPR} to have 3 time scales which change throughout the course of the simulation:
\begin{equation}
  \label{eq:KPR-3scale}
  \begin{pmatrix} u'(t) \\ v'(t) \\ w'(t) \end{pmatrix}
  = \begin{bmatrix} G & e & e \\ e & \alpha & \beta \\ e & -\beta & \alpha \end{bmatrix}
    \begin{pmatrix} \left(u^2-p-2\right)/(2u) \\ \left(v^2-q-2\right)/(2v)\\ \left(w^2-r-2\right)/(2w) \end{pmatrix}
    + \begin{pmatrix} p'(t) / (2u) \\ q'(t) / (2v) \\ r'(t) / (2w) \end{pmatrix}
\end{equation}
over $0< t < 5$, where $p(t) = \frac12 \cos(t)$, $q(t) = \cos(\omega t(1 + e^{-(t-2)^2}))$, and $r(t) = \cos(\omega^2 t(1 + e^{-(t-3)^2}))$.  This problem has analytical solution $u(t) = \sqrt{2+p(t)}$, $v(t) = \sqrt{2+q(t)}$, and $w(t) = \sqrt{2+r(t)}$, stable (but oscillatory) eigenvalues, and its behavior is dictated by the parameters: $e$ that determines the strength of coupling between the time scales, $G<0$ that determines the stiffness at slow time scale, $\alpha$ and $\beta$ that govern oscillations between $v$ and $w$, and $\omega$ that determines the time-scale separation factors between $u$ and $v$, and between $v$ and $w$.
Denoting the full right-hand side vector for \eqref{eq:KPR-3scale} as $\begin{pmatrix} f_u(t,u,v,w) & f_v(t,u,v,w) & f_w(t,u,v,w) \end{pmatrix}^T$, we split the right-hand side into three functions for the slow ($\fslow$), medium ($\fmed$), and fast $\ffast$ scales via
\[
  \underbrace{\begin{pmatrix} f_u(t,u,v,w) \\ 0 \\ 0 \end{pmatrix}}_{\fslow}
  +\underbrace{\begin{pmatrix} 0 \\ f_v(t,u,v,w) \\ 0 \end{pmatrix}}_{\fmed}
  +\underbrace{\begin{pmatrix} 0 \\ 0 \\ f_w(t,u,v,w) \end{pmatrix}}_{\ffast}
\]

Summary statistics from running these nested multirate simulations with the HT-I adaptivity controller at an absolute tolerance $\text{abstol}=10^{-11}$ and a range of relative tolerances $\text{reltol} = \{10^{-2}, 10^{-4}, 10^{-6}, 10^{-8}\}$ are provided in Table \ref{tab:nested-mri}.  All tests used the time scale separation factor $\omega=50$, stiffness factor $G=-10$, and coupling factors $e=5$, $\alpha=-1$, $\beta=1$, the \texttt{ERK22b} MRI method at the slow and intermediate time scales, and ARKODE's default second-order explicit Runge--Kutta method at the fast time scale.  From these, we see that HT-I indeed works well at tracking the dynamics of each time scale, achieving solutions with accuracy metrics within a factor of $\sim\!\!30$ from the requested tolerances, and with work metrics that increase on average by a factor of 33 from one scale to the next.

\begin{table}[!th]
  \centering
  \caption{Summary statistics for multirate simulations of the nested KPR problem \eqref{eq:KPR-3scale} using the HT-I controller at both the slow and intermediate time scales, for various relative tolerances.  The number of slow, intermediate, and fast time steps are shown, as well as the attained accuracy factor as defined in equation \eqref{eq:accuracy}.}
  \begin{tabular}{crrrrr}
    Reltol & Slow Steps & Int Steps & Fast Steps & Accuracy Factor\\
    \hline
    $10^{-2}$ &     84 &      294 &       4,028 & $29.79$ \\
    $10^{-4}$ &  1,081 &   12,896 &     455,531 & $10.19$ \\
    $10^{-6}$ &  1,686 &   91,208 &   2,880,373 & $14.16$ \\
    $10^{-8}$ &  9,793 &  861,182 &  26,943,054 & $6.47$ \\
    \hline
  \end{tabular}
  \label{tab:nested-mri}
\end{table}

\section{Conclusions}
\label{sec:conclusions}

In this work we present two new families of multirate time step adaptivity controllers, \Decoupled\ and \HTol, that are designed to work with embedded MRI methods for adapting time steps when solving problems with multiple time scales.  Comparing against the previously-introduced \Hh\ controllers from \cite{fishAdaptiveTimeStep2023}, the proposed controllers offer dramatically improved performance and flexibility.  While the proposed controller families are able to track multiple time scales for methods of a variety of orders of accuracy, the \Hh\ controllers struggle to select steps with sufficient scale separation between $H$ and $h$, resulting in worse accuracy and higher cost than the proposed families.  The combination of embedded MRI methods with the \Decoupled\ and \HTol\ controllers theoretically support adaptive simulations of problems with an arbitrary number of time scales (here shown for a 3-scale benchmark problem), and achieve high accuracy while maintaining low computational cost.  Of the proposed families, the \HTol\ controllers show much stronger performance for problems where cost is dominated by evaluation of the slow operators, at the expense of requiring estimates of the accumulated fast temporal error.  On the other hand, the \Decoupled\ controllers require no such accumulated fast error estimate, and show the best performance for problems where the fast and slow operators have comparable cost.  Thus, we recommend that readers consider this trade-off when selecting a multirate controller for their application.

We also compared the performance of many adaptive MRI methods on our benchmark problems. Although it is clear that some MRI methods outperform others in these tests, the optimal choice is clearly problem-specific.  We thus encourage practitioners to employ software libraries such as ARKODE \cite{reynolds_arkode_2023} that allow them to easily switch between methods to explore what works best for their application.

A key challenge when using partitioned solvers of any type (e.g., ImEx, MRI) is to construct a good splitting of the ODE right-hand side into components, $f(t,y) = \sum_i f^{\{i\}}(t,y)$.  While we have split the two benchmark problems here to work with explicit, implicit, and ImEx MRI methods, we do not claim that these splittings are optimal for either problem.  Interestingly, although this KPR problem would not generally be considered as stiff, the results in Section \ref{sec:adaptivity_comparisions} indicated a slight preference for implicit and ImEx MRI methods.  Similarly, since the stiff Brusselator benchmark was split such that the stiff terms were subcycled at the fast time scale, it is notable that explicit MRI methods generally achieved the best performance.  This question of how to optimally match explicit, implicit, or ImEx MRI methods to a given application is important but poorly understood.  It is therefore critical for practitioners to additionally experiment with different splittings for their applications.

We note that although this manuscript examined performance on small-scale benchmark problems, the benchmarks were selected to represent some key traits of larger-scale multirate PDE applications.  We are currently testing these proposed MRI adaptivity algorithms on large scale problems arising in simulations of real-time Boltzmann transport equations and tokamak fusion plasmas; however, we leave those results to future publications given the breadth of experiments already performed here.

We further note that not all applications are easily amenable to additive splittings of the form $f(t,y) = \sum_i f^{\{i\}}(t,y)$.  To our knowledge, although some recent work has introduced multirate methods for nonlinearly partitioned systems \cite{buvoli2025multiraterungekuttanonlinearlypartitioned}, those do not yet support either the ``infinitesimal'' structure explored here, or embeddings for temporal error estimation.  However, we anticipate that once those classes of methods have been further developed to provide these features, techniques similar to those in this work can be applied for adapting the step sizes at each time scale.

\corrA{Additionally, the analysis in Section \ref{subsubsec:Hh-analysis} hints at additional criteria that may be used to improve the \Hh\ controller family itself.  Future studies can leverage this criteria to adjust the internal parameters of each \Hh\ controller to ensure their appropriateness for applications with large multirate ratios.}

Finally, we note that the individual single-rate controllers that comprise both the \Decoupled\ and \HTol\ controllers may be chosen independently, and thus different types of single-rate controllers could be used for each component.  For simplicity in this work we consistently selected single-rate component controllers of matching type.  However, future studies may investigate the performance of mixed single-rate controllers, e.g., since it is likely that slow time steps may benefit less from a longer temporal history, a potential combination of an I controller for the slow time scale and a higher-order controller for the fast time scale may be beneficial.

\section*{Acknowledgements}
Daniel R. Reynolds and Vu Thai Luan would like to thank the Vietnam Institute for Advanced Study in Mathematics (VIASM) for their hospitality during the summer research stay, where part of this work was carried out.

\bibliographystyle{elsarticle-num}
\bibliography{main}

\newpage
\appendix

\section{Embedded MERK Methods}
\label{sec:embeddedMERK}

In \cite{luanNewClassHighOrder2020} the authors introduce the methods \texttt{MERK2}, \texttt{MERK3}, \texttt{MERK4}, and \texttt{MERK5}.  Each of these has a beneficial structure that allows an embedding to be included in the next-to-last internal stage, giving rise to the following embedded methods.  As in \cite{luanNewClassHighOrder2020}, we provide only the forcing functions for the internal stages $\left(r_{n,i}(\tau)\right)$, the updated solution $\left(r_n(\tau)\right)$, and the embedding $\left(\tilde{r}_n(\tau)\right)$.  For brevity, we denote the slow function evaluations at each step and stage as $\fslow_{n} = \fslow(t_n,y_n)$ and $\fslow_{n,i} = \fslow(t_n+c_iH, z_i)$, respectively, and we denote the difference functions as $D_{n,i} = \fslow_{n,i}-\fslow_n$.
\begin{itemize}
\item \texttt{MERK21}, that has abscissae $c=\begin{bmatrix} 0 & c_2 & 1 \end{bmatrix}$ and uses the constant $c_2=\tfrac12$:
  \begin{align*}
      r_{n,2}(\tau) &= \fslow_n,\\
      r_n(\tau) &= \fslow_n + \tfrac{\tau}{c_2 H} D_{n,2},\\
      \tilde{r}_n(\tau) &= \fslow_n.
  \end{align*}
\item \texttt{MERK32}, that has abscissae $c=\begin{bmatrix} 0 & c_2 & \tfrac23 & 1 \end{bmatrix}$ and also uses the constant $c_2=\tfrac12$:
  \begin{align*}
      r_{n,2}(\tau) &= \fslow_n,\\
      r_{n,3}(\tau) &= \fslow_n + \tfrac{\tau}{c_2 H} D_{n,2},\\
      r_n(\tau) &= \fslow_n + \tfrac{3\tau}{2 H} D_{n,3},\\
      \tilde{r}_n(\tau) & = r_{n,3}(\tau).
    \end{align*}
\item \texttt{MERK43}, that has abscissae $c=\begin{bmatrix} 0 & c_2 & c_3 & c_4 & c_5 & c_6 & 1 \end{bmatrix}$ and uses the constants $c_2=c_3=\tfrac12$, $c_4=c_6=\tfrac13$ and $c_5=\tfrac56$:
  \begin{align*}
      r_{n,2}(\tau) &= \fslow_n,\\
      r_{n,3}(\tau) &= r_{n,4}(\tau) = \fslow_n + \tfrac{\tau}{c_2 H}D_{n,2},\\
      r_{n,5}(\tau) &= r_{n,6}(\tau) = \fslow_n + \tfrac{\tau}{H}\left( \tfrac{-c_4}{c_3(c_3-c_4)} D_{n,3} + \tfrac{c_3}{c_4(c_3-c_4)} D_{n,4}\right) \\
           &+ \tfrac{\tau^2}{H^2}\left( \tfrac{1}{c_3(c_3-c_4)} D_{n,3} - \tfrac{1}{c_4(c_3-c_4)} D_{n,4} \right),\\
      r_n(\tau) &= \fslow_n + \tfrac{\tau}{H}\left( \tfrac{-c_6}{c_5(c_5-c_6)} D_{n,5} + \tfrac{c_5}{c_6(c_5-c_6)} D_{n,6} \right) \\
           &+ \tfrac{\tau^2}{H^2}\left( \tfrac{1}{c_5(c_5-c_6)} D_{n,5} - \tfrac{1}{c_6(c_5-c_6)} D_{n,6}\right),\\
      \tilde{r}_n(\tau) &= r_{n,5}(\tau).
    \end{align*}
\item \texttt{MERK54}, that has abscissae $c=\left[\; 0\; c_2\; c_3\; c_4\; c_5\; c_6\; c_7\; c_8\; c_9\; c_{10}\; 1\;\right]$ and uses the constants $c_2=c_3=c_5=c_9=\tfrac12$, $c_4=c_6=\tfrac13$, $c_7=\tfrac14$, $c_8=\tfrac{7}{10}$, and $c_{10}=\tfrac23$:
  \begin{align*}
      r_{n,2}(\tau) &= \fslow_n,\\
      r_{n,3}(\tau) &= r_{n,4}(\tau) = \fslow_n + \tfrac{\tau}{c_2 H} D_{n,2},\\
      r_{n,5}(\tau) &= r_{n,6}(\tau) = r_{n,7}(\tau) = \fslow_n + \tfrac{\tau}{H} \left( \alpha_3 D_{n,3} + \alpha_4 D_{n,4}\right)\\
           &+ \tfrac{\tau^2}{H^2}\left( \beta_3 D_{n,3} - \beta_4 D_{n,4} \right),\\
      r_{n,8}(\tau) &= r_{n,9}(\tau) = r_{n,10}(\tau) = \fslow_n + \tfrac{\tau}{H} \left( \alpha_5 D_{n,5} + \alpha_6 D_{n,6} + \alpha_7 D_{n,7}\right)\\
           &- \tfrac{\tau^2}{H^2}\left( \beta_5 D_{n,5} + \beta_6 D_{n,6} + \beta_7 D_{n,7} \right) + \tfrac{\tau^3}{H^3}\left( \gamma_5 D_{n,5} + \gamma_6 D_{n,6} + \gamma_7 D_{n,7} \right),\\
      r_n(\tau) &= \fslow_n + \tfrac{\tau}{H} \left( \alpha_8 D_{n,8} + \alpha_9 D_{n,9} + \alpha_{10} D_{n,10}\right)\\
           &- \tfrac{\tau^2}{H^2}\left( \beta_8 D_{n,8} + \beta_9 D_{n,9} + \beta_{10} D_{n,10} \right)\\
           &+ \tfrac{\tau^3}{H^3}\left( \gamma_8 D_{n,8} + \gamma_9 D_{n,9} + \gamma_{10} D_{n,10} \right),\\
     \tilde{r}_n(\tau) &= r_{n,8}(\tau).
    \end{align*}
    where
    \begin{align*}
      \alpha_3 &= \tfrac{c_4}{c_3 (c_4-c_3)},  \ \alpha_4=\tfrac{c_3}{c_4 (c_3-c_4)},\\
      \alpha_5 &= \tfrac{c_6 c_7}{c_5 (c_5-c_6)(c_5 - c_7)},\  \alpha_6=\tfrac{c_5 c_7}{c_6 (c_6-c_5)(c_6 - c_7)}, \  \alpha_7=\tfrac{c_5 c_6}{c_7 (c_7-c_5)(c_7 - c_6)},\\
      \alpha_8 &= \tfrac{c_9 c_{10}}{c_8 (c_8-c_9)(c_8 - c_{10})},\   \alpha_9=\tfrac{c_8 c_{10}}{c_9 (c_9-c_8)(c_9 - c_{10})},\ \alpha_{10}=\tfrac{c_8 c_{9}}{c_{10} (c_{10}-c_8)(c_{10} - c_{9})}   \\
      \beta_3 &= \tfrac{1}{c_3 (c_3-c_4)},\ \beta_4=\tfrac{1}{c_4 (c_3-c_4)}, \\
      \beta_5 &= \tfrac{c_6+ c_7}{c_5 (c_5-c_6)(c_5 - c_7)},\  \beta_6=\tfrac{c_5 +c_7}{c_6 (c_6-c_5)(c_6 - c_7)}, \  \beta_7=\tfrac{c_5+ c_6}{c_7 (c_7-c_5)(c_7 - c_6)},\\
      \beta_8 &= \tfrac{c_9 +c_{10}}{c_8 (c_8-c_9)(c_8 - c_{10})},\   \beta_9=\tfrac{c_8+ c_{10}}{c_9 (c_9-c_8)(c_9 - c_{10})},\ \beta_{10}=\tfrac{c_8 +c_{9}}{c_{10} (c_{10}-c_8)(c_{10} - c_{9})}\\
      \gamma_5 &= \tfrac{1}{c_5 (c_5-c_6)(c_5 - c_7)},\  \gamma_6=\tfrac{1}{c_6 (c_6-c_5)(c_6 - c_7)}, \  \gamma_7=\tfrac{1}{c_7 (c_7-c_5)(c_7 - c_6)},\\
      \gamma_8 &= \tfrac{1}{c_8 (c_8-c_9)(c_8 - c_{10})},\   \gamma_9=\tfrac{1}{c_9 (c_9-c_8)(c_9 - c_{10})},\ \gamma_{10}=\tfrac{1}{c_{10} (c_{10}-c_8)(c_{10} - c_{9})}.
    \end{align*}
\end{itemize}


\end{document}